\documentclass{amsart}
\usepackage{amsmath,amsthm}
\usepackage{graphics}
\usepackage{color}
\usepackage{epsfig}

\graphicspath{{./Figures/}}

\theoremstyle{plain}

\newtheorem{thm}{Theorem}[section]
\newtheorem{cor}[thm]{Corollary}

\newtheorem{prop}[thm]{Proposition}
\newtheorem{lemma}[thm]{Lemma}
\newtheorem{algor}[thm]{Algorithm}
\newtheorem{claim}[thm]{Claim}
\newtheorem{remark}[thm]{Remark}
\newtheorem{example}[thm]{Example}

\newcommand{\N}{\ensuremath{\mathbb{N}}}
\newcommand{\Q}{\ensuremath{\mathbb{Q}}}

\newcommand{\Z}{\ensuremath{\mathbb{Z}}}
\newcommand{\K}{\ensuremath{{\mathcal K}}}
\newcommand{\D}{\ensuremath{{\mathcal D}}}
\newcommand{\bdry}{\ensuremath{\partial}}
\DeclareMathOperator{\len}{\ell}
\DeclareMathOperator{\SL}{SL}
\DeclareMathOperator{\sgn}{sgn}
\DeclareMathOperator{\Stab}{Stab}
\DeclareMathOperator{\vol}{vol}

\newcommand{\comment}[1]{}
\newcommand{\bit}[1]{\mbox{\boldmath$#1$}}
\newcommand{\mat}[4]{\left( \begin{array}{cc} #1 & #3 \\ #2 & #4 \end{array} \right)}

\theoremstyle{plain}
\newtheorem*{thm_main}{Theorem~\ref{thm:generalapplication}}
\newtheorem*{corA}{Corollary \ref{cor:largevolumemanysurfaces}}
\newtheorem*{corB}{Corollary \ref{cor:largevolumenosurfaces}}

\begin{document}

\title[Surfaces in the complements of large volume Berge knots]{Closed essential surfaces in the complements of large volume Berge knots}

\author{Kenneth L. Baker}
\address{Department of Mathematics, University of Georgia \\ Athens, Georgia 30602}
\email{kb@math.uga.edu}

\thanks{This work was partially supported by a graduate traineeship from the VIGRE Award at the University of Texas at Austin and a VIGRE postdoc under NSF grant number DMS-0089927 to the University of Georgia at Athens.}



\begin{abstract}
We construct an algorithm that lists all closed essential surfaces in the complement of a knot that lies on the fiber of a trefoil or figure eight knot.  Such knots are Berge knots and hence admit lens space surgeries.  Furthermore they may have arbitrarily large hyperbolic volume.  Using this algorithm we concoct large volume Berge knots of two flavors:  those whose complement contains arbitrarily many distinct closed essential surfaces, and those whose complement contains no closed essential surfaces.
  \end{abstract}

\maketitle

\section{Introduction}
Given any knot, it is an interesting problem to determine all closed essential surfaces in its complement.  Understanding such surfaces in the complement of a knot may give insight into properties of the knot itself.  

In \cite{baker:sdavobkI} and \cite{baker:sdavobkII} we show that among the known knots with lens space surgeries, the Berge knots \cite{berge:skwsyls}, only those that lie on the fiber of a trefoil or the figure eight knot may have complements with arbitrarily large hyperbolic volume.  The question of whether the presence of certain essential surfaces in the complements of these large volume Berge knots correlates to this schism motivated this present work.  We answer this question in the negative.  

Let $\K_\eta$ be the collection of knots that lie on the fiber of the left handed trefoil.
\begin{corA}
There exist arbitrarily large volume knots in $\K_\eta$ with arbitrarily many distinct closed essential surfaces in their complements.
\end{corA}
\begin{corB}
There exists arbitrarily large volume knots in $\K_\eta$ with no closed essential surfaces in their complements.
\end{corB}
Similar results hold if $\K_\eta$ is the collection of knots that lie on the fiber of the figure eight knot.

These corollaries come from combining results of \cite{baker:sdavobkI} with the following main result from this paper.
\begin{thm_main}[Short version]
Let $L$ be a knot in $\K_\eta$ whose slope has continued fraction expansion $[b_1, b_2, \dots, b_k]$.  Assume the signs of the coefficients alternate, $|b_i| \geq 2$ for each $i$, and $b_1 \neq 2$.  Then every closed essential surface in the complement of $L$ corresponds to a solution of the following equation:
\[ 0 =
   \sum_{i \in I} -b_i + 
   \sum_{j \in J} b_j +
   \begin{cases}
      0 & \mbox{ if } 1 \in I \\
      -1 & \mbox{ otherwise}
   \end{cases}
\]
where $I$ and $J$ are subsets of $\{1, \dots, k\}$ each not containing consecutive integers and $1 \not \in I \cap J$.  
\end{thm_main}

Furthermore, in Proposition~\ref{prop:differences_in_surfaces} we show that if $\widehat{S}$ is a closed essential surface in the complement of $L \in \K_\eta$, then there is either a longitudinal annulus or a meridional annulus between $\widehat{S}$ and $L$.  If there is a longitudinal annulus, then $\widehat{S}$ bounds a handlebody in $S^3$ containing $L$.  If there is a meridional annulus, then $\widehat{S}$ does not bound a handlebody.

\subsection{Outline of Theorem~\ref{thm:generalapplication}}
The knots of interest here lie on the fiber of a trefoil or the figure eight knot.  Perhaps it is better to conceive of them as lying on a once-punctured torus page of an open book.  By assuming a closed essential surface in the complement of one of these knots has been isotoped to intersect the binding (i.e.\ a trefoil or figure eight knot) minimally, we may exploit the fibration of the complement of the binding to get a grasp on the surface.  Indeed, this fibration is just a once-punctured torus bundle.  

Culler, Jaco, and Rubinstein \cite{cjr:isioptb} as well as Floyd and Hatcher \cite{floydhatcher:isiptb} have developed algorithms that list all properly embedded incompressible surfaces in any given orientable once-punctured torus bundle.  For our purposes we study properly embedded surfaces in once-punctured torus bundles that are disjoint from a non-trivial level curve (i.e.\ an essential simple closed curve on a fiber) and incompressible in the complement of this level curve.  To do so, in Sections \ref{sec:surfacetypes} and \ref{sec:surfacestructure} we adapt \cite{cjr:isioptb} to accommodate the presence of such a level curve.   This leads us to Algorithm~\ref{algor:main_algorithm} which lists all so-called {\em twisted surfaces} in a given once-punctured torus bundle that are essential in the complement of a given level curve.  These twisted surfaces contain all surfaces that serve our ultimate purpose.

Since we are primarily concerned with surfaces in once-punctured torus bundles whose boundaries are meridional curves of bindings of open books, in Section~\ref{sec:framing} we retool the discussion on framings of once-punctured torus bundles in \cite{cjr:isioptb} to more readily detect when these twisted surfaces have meridional boundary components.  Thereby Algorithm~\ref{algor:framing_algorithm} extends Algorithm~\ref{algor:main_algorithm} to list only those surfaces with the desired meridional boundaries.

In Section~\ref{sec:closedsurfaces} we show that the same conditions which imply twisted surfaces with meridional boundaries are essential in the complement of level curves ensure that the corresponding capped off surfaces in the open book remain incompressible in the complement of the level curve when the binding is filled in.  In Section~\ref{sec:surfacecomplements} we examine the complements of these closed surfaces and their relationship to the knot.

In Section~\ref{sec:bergeknotalg} Algorithm~\ref{algor:streamlined_alg} streamlines these algorithms for the Berge knots that lie on the fiber of a trefoil or the figure eight knot.  We then give examples of direct applications of these algorithms.  This streamlined algorithm and its following examples exhibit the process by which the equation of Theorem~\ref{thm:generalapplication} is obtained.  

\subsection{Acknowledgements}
The author wishes to thank John Luecke for his direction an many useful conversations.

\section{Preliminaries}

A simple closed curve on a compact surface is {\em essential} if on the surface it neither bounds a disk nor is isotopic into the boundary.  A properly embedded surface in a compact $3$-manifold is {\em essential} if it is incompressible, $\bdry$--incompressible, and not $\bdry$--parallel.  A compact $3$-manifold that contains no closed essential surfaces is called {\em small}.

\subsection{The once-punctured torus, $\SL_2(\Z)$, and bundles}\label{subsec:opt}
We recall some notation and terminology of \cite{cjr:isioptb}.

\begin{figure}
\centering
\input{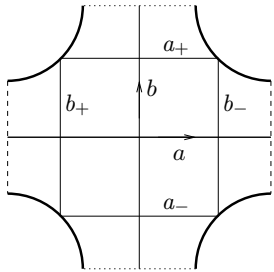}
\caption{The once-punctured torus $T$.}
\label{fig:o-pt}
\end{figure}

Let $T$ be the once-punctured torus shown (with the appropriate identifications of edges) in Figure~\ref{fig:o-pt}.  Also shown are the oriented simple closed curves $a$ and $b$ and arcs $a_+$, $a_-$, $b_+$, and $b_-$.  Orient $T$ so that $a \cdot b = +1$.

Let $\alpha$ and $\beta$ be the left-handed Dehn twists along $a$ and $b$ respectively. Then $\alpha$ and $\beta$ generate the group $\mathcal{H}^+(T)$ of orientation preserving homeomorphisms up to isotopy.  This group is naturally identified with the group of orientation preserving homomorphisms of $H_1(T)$, namely $\SL_2(\Z)$.  Assuming $H_1(T)$ is generated by $[a]$ and $[b]$, this identification is done by 
\[
\alpha \mapsto A = \begin{pmatrix} 1 & -1 \\ 0 & 1 \end{pmatrix} \mbox{ and } 
\beta \mapsto B = \begin{pmatrix} 1 & 0 \\ 1 & 1 \end{pmatrix}.
\]
These two matrices generate $\SL_2(\Z)$.
The two matrices 
\[ P = \begin{pmatrix} 0 & -1 \\ 1 & 0 \end{pmatrix} = BAB \mbox{ and } 
Q = \begin{pmatrix} 1 & 1 \\ -1 & 0 \end{pmatrix} = B^{-1}A^{-1}\]
also generate $\SL_2(\Z)$.
This gives the presentation
\[ \SL_2(\Z) = \langle P, Q | P^2=Q^3=-I \rangle.\]

Note that $P$ is identified with the isotopy class of the homeomorphism $\phi$ of $T$ that is induced by rotating Figure~\ref{fig:o-pt} by $\frac{\pi}{2}$.

If $\eta$ is an orientation preserving homeomorphism of $T$, then we may form the once-punctured torus bundle 
\[ M = T \times I /_\eta \]
defined by identifying $T \times \{1\}$ with $T \times \{0\}$ of $T \times I$ by $(x,1) \mapsto (\eta(x), 0)$.

If we express $\eta$ as a composition of such homeomorphisms $\eta = \eta_n \eta_{n-1} \dots \eta_2 \eta_1$, then $M$ may be divided into blocks as
\[ M = T \times I/_{\eta_1} T \times I/_{\eta_2} \dots /_{\eta_{n-1}} T \times I/_{\eta_n} \]
where $T \times \{1\}$ of the $i$th block is identified with $T \times \{0\}$ of the $i+1$th (mod~$n$) block according to $(x,1)_i \mapsto (\eta_i(x), 0)_{i+1}$.

If $H_i \in \SL_2(\Z)$ is identified with the isotopy class of the homeomorphism $\eta_i$, then $M$ has characteristic class $[H_n H_{n-1} \dots H_2 H_1]$.

\subsection{Continued fractions}

We define a {\em continued fraction expansion}  $[b_1, b_2, \dots, b_k]$ of $\frac{p}{q}$ as follows:
\[ \frac{p}{q} = [b_1, b_2, \dots, b_k] = 
\cfrac{1}{b_1-
 \cfrac{1}{b_2-
  \cfrac{1}{\dots -
   \cfrac{1}{b_k
}}}}
\]  
where $b_i \in \Z$ for $i = 1, \dots, k$ are the {\em coefficients}.  For each $i = 1, \dots, k$, the {\em $i$th partial fraction} of $[b_1, b_2, \dots, b_k]$ is $[b_1, b_2, \dots, b_i]$.
The $0$th partial fraction the continued fraction $[\quad]$ which represents $0$.  

A rational number has many continued fraction expansions.  
We say $[b_1, b_2, \dots, b_k]$ is a {\em minimal continued fraction expansion (MCFE)} for $\frac{p}{q}$ if $|b_i| \geq 2$ for all $i = 1, \dots, n$. 
In order for $\frac{p}{q}$ to have a MCFE, it must be that $|\frac{p}{q}| < 1$.
We say $[b_1, b_2, \dots, b_k]$ is the {\em simple continued fraction expansion (SCFE)} for $\frac{p}{q}$ if the coefficients alternate signs, $|b_k| \geq 2$, and $b_i \neq 0$ for $i \geq 2$.
Though continued fraction expansions are not unique, the Euclidean algorithm gives each rational number $\frac{p}{q}$ a unique SCFE.  Furthermore, the SCFE for $\frac{p}{q}$ has $0$ as its first coefficient if and only if $|\frac{p}{q}| \geq 1$ (see Lemma~\ref{lem:scfesize}).

Let $\mathcal{D}$ be the {\em Farey Diagram} shown in Figure~\ref{fig:diagram}(a).  $\mathcal{D}$ is a disk with the extended rational numbers marked on its boundary.  An edge joins vertices $\frac{a}{b}$ to $\frac{c}{d}$ if and only if $ad-bc = \pm 1$.  The edge from $\frac{a}{b}$ to $\frac{c}{d}$ is the ``long'' edge of the triangle whose third vertex is $\frac{a+c}{b+d}$.  

\begin{figure}
\centering
$\begin{array}{ccc}
\input{Figures/Diagram.pstex_t} & \mbox{\quad}&\input{Figures/Diagram-Example.pstex_t}\\
\mbox{(a) The diagram $\mathcal{D}$.} && \mbox{(b) The edge-path $E_{[2,-1]}$.}
\end{array}$
\caption{}
\label{fig:diagram}
\end{figure}

We utilize the connections between $\mathcal{D}$, $\SL_2(\Z)$, and continued fraction expansions in the vein of Floyd-Hatcher \cite{floydhatcher:isiptb} and Hatcher-Thurston \cite{hatcherthurston:isi2bkc} and borrow some of their terminology.  One may also look to Kirby-Melvin \cite{kirbymelvin:dsmiatsc} for a more thorough treatment of these connections.

Let $[\bit{b}]$ denote the continued fraction $[b_1, \dots, b_k]$.  An {\em (oriented) edge-path}  $E_{\bit{b}} = \{e_0, e_1, e_2, \dots, e_k\}$ from $\frac{1}{0}$ to $\frac{p}{q}$ corresponds uniquely to a continued fraction expansion $[b_1, b_2, \dots, b_k]$ for $\frac{p}{q}$ where $e_0$ is the edge from $\frac{1}{0}$ to $\frac{0}{1}$ and the end points of each edge $e_i$ are the partial fractions $\frac{p_{i-1}}{q_{i-1}} = [b_1, \dots, b_{i-1}]$ and $\frac{p_i}{q_i}=[b_1, \dots, b_i]$.  At the vertex $\frac{p_{i-1}}{q_{i-1}}$, facing inward, the edge $e_i$ is $|b_i|$ triangles apart from $e_{i-1}$, to the left if $b_i < 0$ and to the right if $b_i > 0$.  If $b_i=0$ the edge $e_i$ retraces the edge $e_{i-1}$.
Figure~\ref{fig:diagram}(b) shows an example of the edge-path corresponding to $[2,-1]$.

An edge-path $E_{\bit{m}}$ associated to a MCFE $[\bit{m}]$ is a {\em minimal edge-path}.  This edge-path is minimal in the sense that no triangle of $\mathcal{D}$ has two consecutive edges of $E_{\bit{m}}$ on its boundary.

\begin{figure}
\centering
\input{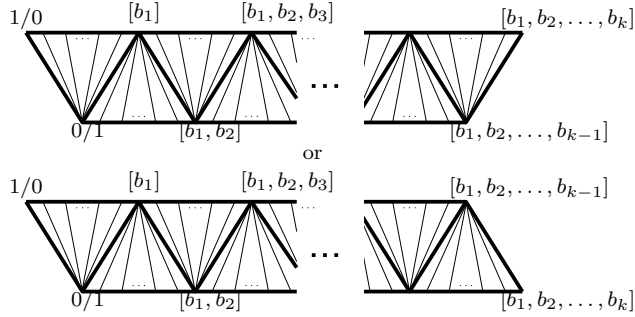}
\caption{The subcomplex of $\mathcal{D}$ associated to a SCFE.}
\label{fig:simple_edge-path}
\end{figure}

Minimal edge paths from $\frac{1}{0}$ to $\frac{p}{q}$ are contained in the finite subcomplex $E_{\frac{p}{q}}$ of $\mathcal{D}$ comprised of the edge-path associated to the SCFE $[a_1, a_2, \dots, a_l]$ for $\frac{p}{q}$ and the triangles between its edges.  See Figure~\ref{fig:simple_edge-path}.  Observe that a minimal edge path may only involve the heavier edges of the larger triangles.  As noted in \cite{hatcherthurston:isi2bkc}, the number of minimal edge paths from $\frac{1}{0}$ to $\frac{p}{q}$ (and hence the number of MCFEs for $\frac{p}{q}$) is finite.

The following two lemmas will be used in \S\ref{sec:generalities}
\begin{lemma}\label{lem:scfesize}
Let $[\bit{b}] = [b_1, b_2, \dots, b_k]$ be the SCFE for $\frac{p}{q} \neq \pm1$.  If $b_1 \neq 0$ then $|\frac{p}{q}| < 1$ and $\sgn(\frac{p}{q}) = \sgn(b_1)$.  If $b_1 = 0$ then $|\frac{p}{q}| > 1$ and $\sgn(\frac{p}{q}) = -\sgn(b_2)$.
\end{lemma}
\begin{proof}
This is immediate when one considers the edge-path in $\mathcal{D}$ associated to a SCFE.
\end{proof}

\begin{lemma}\label{lem:threesums}
If the three final partial fractions $p_l/q_l$, $p_{l-1}/q_{l-1}$, and $p_{l-2}/q_{l-2}$ of $[a_1, \dots, a_{l-2}, a_{l-1}, a_l]$ equal the three final partial fractions $r_k/s_k$, $r_{k-1}/s_{k-1}$, and $r_{k-2}/s_{k-2}$ of $[b_1, \dots, b_{k-2}, b_{k-1}, b_k]$ respectively, then $a_l = b_k$.
\end{lemma}
\begin{proof}
The three final partial sums determine the final two edges of the associated edge-paths.  With these partial sums equal, the two edges are the same, and hence they have the same number of triangles between them.  This number is the final coefficient $a_l=b_k$.
\end{proof}

\subsubsection{Relating continued fraction expansions}\label{subsubsec:rcfe}

Let $[\bit{a}] = [a_1, a_2, \dots, a_l]$ and $[\bit{b}] = [b_1, b_2, \dots, b_k]$ be two continued fraction expansions for the same rational number.  Then $[\bit{a}]$ may be obtained from $[\bit{b}]$ by a finite sequence of the following {\em elementary moves} and their inverses:
\begin{align*}
\tag{CF1} [\dots, r, s, \dots] &\mapsto [\dots, r \pm 1, \pm 1, s \pm 1, \dots], \\
\tag{CF1'}[\dots, r] &\mapsto [\dots, r \pm 1, \pm 1], \\
\tag{CF2}[\dots, r+s, \dots] &\mapsto [\dots, r, 0, s, \dots], \\
\tag{CF2'}[\dots, r] &\mapsto [\dots, r, s, 0]. \\
\end{align*}
These may be understood via the corresponding moves on edge-paths $E_{\bit{a}}$ and $E_{\bit{b}}$.  Note that (CF2') may be obtained from the other three moves.

Let $[\bit{a}] = [a_1, a_2, \dots, a_k]$ be the SCFE for $\frac{p}{q}$.  Then any MCFE for $\frac{p}{q}$ may be obtained from $[\bit{a}]$ by a sequence of moves on non-adjacent coefficients $a_i \neq 0$ of $[\bit{a}]$ of the following forms
\begin{align*}
\tag{M} [\dots, a_{i-1}, a_i, a_{i+1}, \dots] &\mapsto 
[\dots, a_{i-1} \pm 1, \underbrace{\pm 2, \pm 2, \dots, \pm 2}_{ |a_i - 1| }, a_{i+1} \pm 1, \dots], 
\end{align*}
for $i \neq 1$ or $k$ and
\begin{align*}
\tag{M'} [\dots, a_{k-1}, a_k] & \mapsto [\dots, a_{k-1} \pm 1, \underbrace{\pm 2, \pm 2, \dots, \pm 2}_{ |a_k - 1| }].
\end{align*}
When (M) or (M') is applied to a coefficient $a_i = \pm 1$, it is tantamount to applying the inverses of the move (CF1) or (CF1') respectively.  We use the ``$+$'' if $a_i < 0$ and the ``$-$'' if $a_i > 0$.  The moves (M) and (M') may be obtained as repeated applications of move (CF1) together with one application of move (CF1') for (M').  Again, these may be understood via the corresponding edge-paths.

\begin{lemma}\label{lem:matrix_scfemcfe}
Assume $[\bit{b}] = [b_1, \dots, b_k]$ is a SCFE for $\frac{p}{q}$ and $[\bit{a}] = [a_1, \dots, a_l]$ is a MCFE for $\frac{p}{q}$.  If $\frac{r}{s} = [\bit{b},N] = [\bit{a}, N']$ for integers $N$ and $N'$, then either $N' = N$ or $N' = N-\sgn(b_k)$.
\end{lemma}

\begin{proof}
  Consider the edge-paths in $\D$ corresponding to these two continued fraction expansions of $\frac{r}{s}$.  Since both of these continued fractions have the same penultimate partial sum of $\frac{p}{q}$, their edge-paths $E_{\bit{b},N}$ and $E_{\bit{a},N'}$ share the same final edge.  The integers $N$ and $N'$ describe how far apart are the final edges of $E_{\bit{b},N}$ and $E_{\bit{a},N'}$ from the final edges of $E_{\bit{b}}$ and $E_{\bit{a}}$ respectively.  

Because any MCFE for a rational number is obtained from its SCFE by a sequence of moves (M) and (M'), the final edges of $E_{\bit{b}}$ and $E_{\bit{a}}$ either coincide or are distinct edges of a triangle in $\D$.
If these edges coincide, then $N=N'$.  If the edges are distinct, then as apparent from Figure~\ref{fig:simple_edge-path} the final edge of $E_{\bit{a}}$ is a $-\sgn(b_k)$ triangle apart from the final edge of $E_{\bit{b}}$.  Thus $N = N'-\sgn(b_k)$.
 \end{proof}



\subsubsection{Coefficient sums and lengths}

If $[\bit{b}] = [b_1, b_2, \dots, b_k]$ is a continued fraction expansion, then  $\sigma(\bit{b}) = \sum_{i=1}^k b_i$ is the {\em coefficient sum} of $[\bit{b}]$.
 
 \begin{lemma}\label{lem:coeffsums}
Assume $[\bit{b'}]$ is obtained from $[\bit{b}]$ by an elementary move.
Then the coefficient sums differ as follows:
\begin{align*}
\sigma(\bit{b'}) - \sigma(\bit{b}) = 
\begin{cases} 
\pm 3 \mbox{ if } [\bit{b'}] \mbox{ is obtained by (CF1)}, \\
\pm 2 \mbox{ if } [\bit{b'}] \mbox{ is obtained by (CF1')},\\ 
0     \mbox{ if } [\bit{b'}] \mbox{ is obtained by (CF2)}. \\ 
\end{cases} 
\end{align*}
The difference varies if $[\bit{b'}]$ is obtained from $[\bit{b}]$ by move (CF2').

If $[\bit{b'}]$ is obtained from $[\bit{b}]$ by move (M) at $b_i$ for $1<i<k$ or move (M') at $b_k$, then
\[ 
\sigma(\bit{b'}) - \sigma(\bit{b}) =
\begin{cases}
-3 b_i &\mbox{ if } i \neq k, \\
-3 b_k +1 &\mbox{ if } b_k >0, \\
-3 b_k -1 &\mbox{ if } b_k <0. 
\end{cases}
\]
\end{lemma}
\begin{proof}
The proof is immediate.
\end{proof}

If $[\bit{b}] = [b_1, b_2, \dots, b_k]$ is a continued fraction expansion, then  $\len(\bit{b}) = k$ is the {\em length} of $[\bit{b}]$.

\begin{lemma}\label{lem:length}
Assume $[\bit{b'}]$ is obtained from $[\bit{b}]$ by an elementary move.
Then the lengths differ as follows:
\begin{align*}
\len(\bit{b'}) - \len(\bit{b}) = 
\begin{cases} 
1 \mbox{ if } [\bit{b'}] \mbox{ is obtained by (CF1) or (CF1')}, \\
2     \mbox{ if } [\bit{b'}] \mbox{ is obtained by (CF2) or (CF2')}. \\ 
\end{cases} 
\end{align*}
Also if $[\bit{b'}]$ is obtained from $[\bit{b}]$ by move (M) or (M') at $b_i$ then
\[ 
\len(\bit{b'}) - \len(\bit{b}) = |b_i| - 2.
\]
\end{lemma}
\begin{proof}
The proof is immediate.
\end{proof}



\subsection{$\SL_2(\Z)$ and continued fraction expansions}
 

\begin{lemma} \label{lem:matrix_cfe1}
Let $W =\mat{x}{y}{t}{u} \in \SL_2(\Z)$.  Then
\begin{enumerate}
\item 
\[
\begin{array}{ccc}
W = \pm B^{n_1} A^{n_2} \dots B^{n_k} & & \frac{x}{y} = [n_1, n_2, \dots, n_k] \\
\mbox{or} & \mbox{if and only if}&\mbox{and} \\ 
W=\pm B A B A^{n_1} B^{n_2} \dots B^{n_k} & &\frac{t}{u} = [n_1, n_2, \dots, n_{k-1}], \\
\end{array} 
\]
and
\item 
\[
\begin{array}{ccc}
W = \pm B^{n_1} A^{n_2} \dots A^{n_k} & & \frac{x}{y} = [n_1, n_2, \dots, n_{k-1}] \\
\mbox{or} & \mbox{if and only if}&\mbox{and} \\ 
W=\pm B A B A^{n_1} B^{n_2} \dots A^{n_k} & & \frac{t}{u} = [n_1, n_2, \dots, n_{k}]. \\
\end{array} 
\]
\end{enumerate}
\end{lemma}

\begin{proof}
The proof is straightforward.
\end{proof}

\begin{lemma}~\label{lem:matrix_cfe2}
Let $W = \mat{x}{y}{t}{u} \in \SL_2(\Z)$ and $\frac{x}{y} = [n_1, n_2, \dots, n_k]$.  Then depending on the parity of $k$, 
\[W = \pm B^{n_1} A^{n_2} \dots B^{n_k} A^N \mbox{ or }\pm B A B A^{n_1} B^{n_2} \dots B^{n_k} A^N\]
for some $N \in \Z$.
\end{lemma}

\begin{proof}
Let $W_0 = \mat{x}{y}{t_0}{u_0}$ where $\frac{t_0}{u_0}$ is in lowest terms, has continued fraction expansion $[n_1, n_2, \dots, n_{k-1}]$ and the signs of $t_0$ and $u_0$ are chosen so that $W_0 \in \SL_2(\Z)$.  Then 
\begin{align*}
 W_0^{-1} W &= \mat{u_0}{-y}{-t_0}{x} \mat{x}{y}{t}{u} \\
            &= \mat{1}{0}{t u_0 -u t_0}{1} = A^{u t_0 - t u_0}.
\end{align*}
Thus $W = W_0 A^N$ where $N = u t_0 - t u_0$.  Applying Lemma~\ref{lem:matrix_cfe1} to $W_0$ finishes the lemma.
\end{proof}

Assume $W \in \SL_2(\Z)$ may be expressed as 
\[W = P^J C^{n_k} \dots B^{n_2} A^{n_1}\]
where $|n_i| \geq 2$ for $i = 2, \dots, k$ such that if $k$ is odd then $J \in \{0, 2\}$ and $C = B$, and if $k$ is even then $J \in \{-1, +1\}$ and $C = A$.  We call such an expression for $W$ a {\em special form}.  Contrast this with the special forms as defined in \cite{cjr:isioptb}.

\begin{lemma}\label{lem:specialformMCFE}
$W =\mat{x}{y}{t}{u}\in \SL_2(\Z)$ has a special form $P^J C^{n_k} \dots B^{n_2} A^{n_1}$ if and only if $\frac{x}{y}$ has MCFE $[n_k, \dots, n_2]$.
\end{lemma}
\begin{proof}
This is a consequence of Lemmas~\ref{lem:matrix_cfe1} and \ref{lem:matrix_cfe2}.
\end{proof}

\subsection{Curves on the once-punctured torus}

Let $\K$ be the set of isotopy classes of unoriented essential simple closed curves on the once-punctured torus $T$.  If $c \in \K$ then, after choosing an orientation for $c$, $[c] = p [a] + q [b] \in H_1(T)$ for relatively prime integers $p$ and $q$.  Since $c$ with the opposite orientation has homology $-p[a]+-q[b]$, we may then identify the members $c \in \K$ with $\Q \cup \infty$ via the correspondence $c \mapsto p/q$.  These isotopy classes and their corresponding rational numbers are called {\em slopes}.

\begin{lemma}\label{lem:KandCFE}
Let $K \in \K$ have slope $\frac{p}{q}$.
Then $\frac{p}{q}$ has a continued fraction expansion $[r_n, \dots, r_2, r_1]$ of odd length if and only if
\[ 
K = \beta^{r_n} \circ \cdots \circ \alpha^{r_2} \circ \beta^{r_1} (a).
\]
\end{lemma}

\begin{proof}
Assume $[r_n, \dots, r_2, r_1]$ is a continued fraction expansion for $\frac{p}{q}$ of odd length.  
Let $W \in \SL_2(\Z)$ be the change of basis matrix $W = \mat{p}{q}{p'}{q'}$ where $\frac{p'}{q'} = [r_n, \dots, r_2]$.  Since $n$ is odd, by Lemma~\ref{lem:matrix_cfe1}, $W = \pm B^{r_n} \dots A^{r_2} B^{r_1}$

Via the correspondence between $H_1(T)$ and homeomorphisms of $T$, we have that up to an orientation on $K$
\[ K = \beta^{r_n} \circ \dots \circ \alpha^{r_2} \circ \beta^{r_1} \circ \alpha^N (a), \]
  since $[K] = W [a]$.  Because $\alpha^N (a) = a$, 
\[ K = \beta^{r_n} \circ \dots \circ \alpha^{r_2} \circ \beta^{r_1} (a). \]

These steps all reverse for the other implication. 
\end{proof}

\begin{remark}\label{rem:KandCFE}
If 
\[K = \alpha^{r_n} \circ \cdots \circ \alpha^{r_2} \circ \beta^{r_1} (a), \]
then 
\[K = \beta^0 \circ \alpha^{r_n} \circ \cdots \circ \alpha^{r_2} \circ \beta^{r_1} (a). \]
Hence the slope of $K$ has the continued fraction expansion $[0, r_n, \dots, r_2, r_1]$.
\end{remark}

\section{Twisted surfaces}\label{sec:surfacetypes}

We begin by describing a certain type of surface in once-punctured torus bundles that may be disjoint from a level knot.  Much of the terminology and methods used here are borrowed or adapted from \cite{cjr:isioptb}.

 \subsection{Construction of twisted surfaces}

Here we follow parts of \cite[\S 2]{cjr:isioptb}.  We describe the so-called {\em twisted surfaces} in a once punctured torus bundle that are disjoint from a given essential level knot.  We then prove when they are essential.

\begin{figure}%
\centering
\input{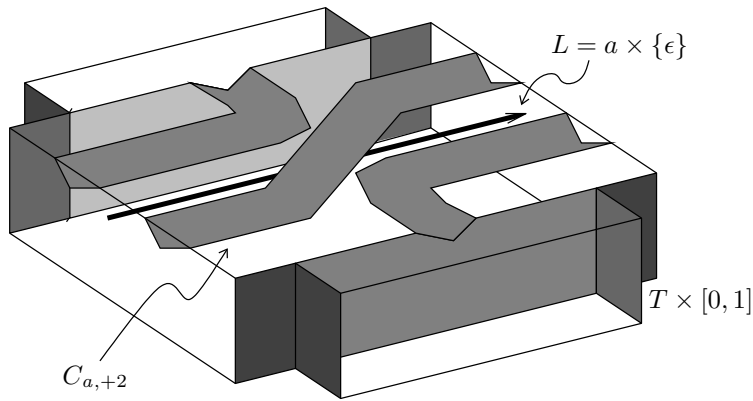}
\caption{The twisted saddle $C_{{\mathit a},+2}$ with level knot $a \times \{\epsilon\}$.}
\label{fig:Ca+2_with_knot}
\end{figure}

To describe the entire gamut of essential surfaces in once-punctured torus bundles, Culler-Jaco-Rubinstein begin by listing several kinds of surfaces embedded in $T \times [0, 1]$.  For our purposes we need only to discuss the {\em twisted saddle} surfaces $C_{{\mathit a},n}$ and $C_{{\mathit b},n}$ (cf.\ \cite[\S 2.1]{cjr:isioptb}).  The surface $C_{{\mathit a},n}$ embeds in $T \times [0, 1] - N(a \times \{\epsilon\})$ where $\epsilon > 0$ (so that $a \times \{\epsilon\}$ is a slight upward push-off of $a \times \{0\}$).  See Figure~\ref{fig:Ca+2_with_knot}.

As in \cite[\S 2.5]{cjr:isioptb}, let $n(1), n(2), \dots, n(k)$ be integers and $J \in \{-1,0 ,1 ,2\}$ have the same parity as $k$.  Consider the bundle
\begin{align*}
M &= T \times I /_{\alpha^{n(1)}} T \times I /_{\beta^{n(2)}} \dots T \times I/_{\gamma^{n(k)}} T\times I /_{\phi^J} \\
 &= T \times I /_{\phi^J \gamma^{n(k)}  \dots \beta^{n(2)} \alpha^{n(1)}}
\end{align*}
where $\gamma = \beta$ if $k$ is even and $\gamma = \alpha$ if $k$ is odd.
As written, $M$ contains $k+1$ blocks, but we will often consider the $k+1$st block as part of the $k$th block.  Let $L = a \times \{\epsilon\}$ in the first block of $M$ for small $\epsilon > 0$ be the level knot at hand.  We may actually think of $L$ as being $a \times \{0\}$, but it is useful to have $L$ not on a fiber along which blocks of $M$ are glued.   

Construct the surface $R$ by putting twisted saddles in the first $k$ blocks of $M$, $C_{{\mathit a},n(i)}$ if $i$ is odd and $C_{{\mathit b},n(i)}$ if $i$ is even, and the vertical disks ${\mathit a}_\pm \times I$ or ${\mathit b}_\pm \times I$ in the $(k+1)$th block of $M$.  These surfaces fit together to make a properly embedded connected surface $R$ which is disjoint from $L$.  $R$ is orientable if $k$ is even and non-orientable if $k$ is odd.

\begin{figure}
\centering
\input{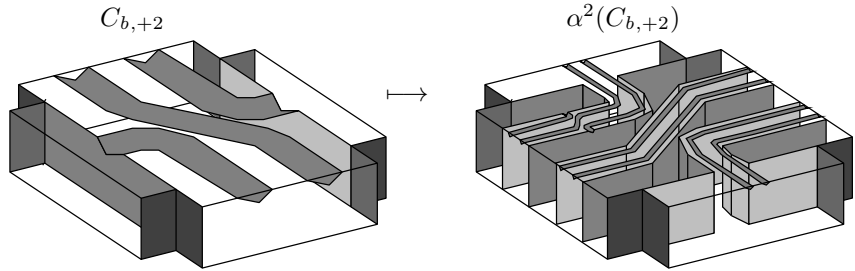}
\caption{``Pushing'' $\alpha^2$ through $C_{{\mathit b},+2}$}
\label{fig:twistCb}
\end{figure}

\begin{figure}
\centering
\input{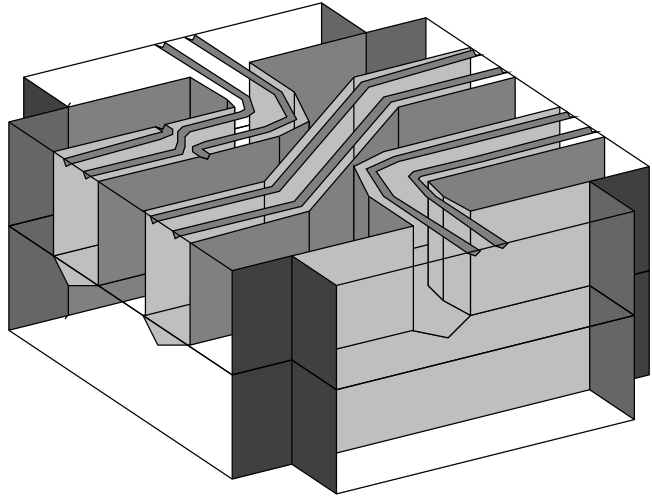}
\caption{The twisted saddles $C_{{\mathit a},+2}$ and $C_{{\mathit b},+2}$ joined in $T \times I/_{\alpha^2} T \times I$}
\label{fig:joined}
\end{figure}

For an example of how these twisted saddles fit together to give a properly embedded connected surface, consider 
\[T \times I/_{\alpha^2} T \times I.\]
View this as two blocks of a once-punctured torus bundle joined together by the homeomorphism $\alpha^2$.  Let the first block (the lower one) contain a copy of $C_{{\mathit a},+2}$ (as in Figure~\ref{fig:Ca+2_with_knot} without $L$) and the second block contain a copy of $C_{{\mathit b},+2}$.  To attach, we ``push'' the homeomorphism $\alpha^2$ through the second block as in Figure~\ref{fig:twistCb}.
Since $\alpha^2(C_{{\mathit b},+2}) \cap (T \times \{0\})$ now coincides with $C_{{\mathit a},+2} \cap (T \times \{1\})$, 
we may join the blocks together to get Figure~\ref{fig:joined}.

As in \cite[\S 2.5]{cjr:isioptb}, define the {\em twisted surface}  $C(J;n(k), \dots, n(1))$ to be $R$ if $k$ is even and $\bdry N(R)$ if $k$ is odd.
As noted in \cite{cjr:isioptb}:
\begin{itemize}
\item  the surface $C(0;n(k), \dots, n(1))$ has genus $\frac{1}{2} k - 1$ and four boundary components each of which intersects the fiber exactly once,
\item  the surface $C(2;n(k), \dots, n(1))$ has genus $\frac{1}{2} k$ and two boundary components each of which intersects the fiber exactly twice, and
\item  the surfaces $C(-1;n(k), \dots, n(1))$ and $C(+1;n(k), \dots, n(1))$ both have genus $k$ and two boundary components each of which intersects the fiber exactly four times.
\end{itemize}

 Let $S$ and $S'$ be properly embedded connected orientable surfaces in once-punctured torus bundles $M$ and $M'$ respectively that are disjoint from essential level knots $L$ and $L'$ respectively.  Then we say $S$ and $S'$ are of the {\em same type} if there is a bundle equivalence from $M$ to $M'$ that maps $S$ to $S'$ and $L$ to $L'$ (cf.\ \cite[\S 2.6]{cjr:isioptb}.)
 If the bundle $M$ contains a surface of type $C(J;n(k), \dots, n(1))$, then $M$ has the characteristic class 
\[ [P^J B^{n(k)} \dots B^{n(2)} A^{n(1)} ].\]

 \subsection{Essential twisted surfaces}

\begin{prop}
[cf.\ Proposition~2.5.1, \cite{cjr:isioptb}]
\label{prop:2.5.1}
Let $J$, $n(1), n(2), \dots, n(k)$, $L$, and $M$ be as above.  The surface $C(J; n(k), \dots, n(2), n(1))$ is essential in $M-N(L)$ if and only if $|n(i)| \geq 2$ for $i = 2, \dots, k$.
\end{prop}

\begin{proof}

We cite the proof of \cite[Proposition 2.5.1]{cjr:isioptb} and show only the parts where we must diverge.  

Let $\widetilde{M}$ be the cyclic cover, corresponding to the fiber, of the bundle
\[ M = T \times I /_{\alpha^{n(1)}} T \times I /_{\beta^{n(2)}} \dots T \times I /_{\varphi^J}. \]
Let $\widetilde{L}$ be the inverse image of $L$ under the covering projection and 
\[S \subset \widetilde{M} - N(\widetilde{L}) \subset \widetilde{M}\] 
be a component of the inverse image of 
\[C(J;n(k),\dots,n(1)) \subset M-N(L) \subset M\]
under the covering projection.  As noted in \cite{cjr:isioptb}, it suffices to show that $S$ is incompressible in ${\widetilde{M}-N(\widetilde{L})}$.  Furthermore $\widetilde{M}$ is divided into blocks which are inverse images of the blocks of $M$, and each block in $\widetilde{M}$ meets $S$ in one disk as in \cite[Figure 4]{cjr:isioptb}.  Let $F$ be the union of the fibers along which the blocks in $\widetilde{M}$ are joined.

Consider the family of all compressing and boundary compressing disks for $S$ in $\widetilde{M}-N(\widetilde{L})$.  This is equivalent to considering the family of compressing disks for $S$ in $\widetilde{M}$ that are disjoint from $\widetilde{L}$.  
By \cite[Proposition 2.5.1]{cjr:isioptb}, the family of all compressing and boundary compressing disks for $S$ in $\widetilde{M}$ (disregarding $\widetilde{L}$) is non-empty if and only if $|n(i)| < 2$ for some $i = 1, \dots, k$.   


Among the disks in this family, consider one that has minimal intersection with $F$.
From the proof of Proposition 2.5.1, such a disk is contained in the solid torus formed by cutting two adjacent blocks of $\widetilde{M}$ along $S$ and joining two of the resulting solid torus components along the annulus of their common intersection with $F$.  If the two blocks are preimages of the $(i-1)$th and $i$th (modulo $k$) blocks of $M$ under the covering projection, then $|n(i)| < 2$ if and only if such a disk exists in this solid torus.  The disk is isotopic to a meridional disk of the solid torus and non-trivially intersects the core of the solid torus.  If $n(i)=0$, then the disk is a compressing disk.  If $n(i) = \pm 1$, then the disk is a boundary compressing disk.  In the case that $i=1$, a component of $\widetilde{L}$ is isotopic to the core of this solid torus (since it is isotopic to the core of the gluing annulus). Thus it will non-trivially intersect the disk.
\end{proof}

 \subsection{Classification of essential twisted surfaces}

A given once-punctured torus bundle with an essential level knot may contain several surfaces essential in the complement of the level knot which all have the type of a twisted surface.  We may determine when two are in the same isotopy class.  This is effectively done in \cite[\S 4.1]{cjr:isioptb}.

Recall that if $M$ is a fiber bundle, then an isotopy of $M$ which is a bundle equivalence at each time is called a {\em bundle isotopy}.

\begin{prop} \label{prop:uniqueness}
Let $M$ be a once-punctured torus bundle containing surfaces $S$ and $S'$ of types $C(J; n(k), \dots, n(1))$ and $C(J'; m(k'), \dots, m(1))$, respectively, that are disjoint from an essential level curve $L$.  Then $S$ and $S'$ are isotopic in the complement of $L$ if and only if $k = k'$, $J = J'$, and $(m(k), \dots, m(1)) = (n(k), \dots, n(1))$.   Moreover there is a bundle isotopy taking $(S, L)$ to $(S', L)$.
\end{prop}

\begin{proof}
This follows directly from the proof of \cite[Proposition~4.1.3]{cjr:isioptb}.  Note that because of $L$, we have that $(m(k), \dots, m(1))$ and $(n(k), \dots, n(1))$ are related by equality rather than cyclic permutation. 
\end{proof}

We now relate essential twisted surfaces in a based once-punctured torus bundle to certain expressions of the bundle's characteristic class.

\begin{prop}\label{prop:existence}
Let $M = T \times I /_\eta$ be a once-punctured torus bundle containing an essential level curve $L \subset T \times \{0\}$.
Assume $\eta$ is identified with $H \in \SL_2(\Z)$.  Then for each element $W \in \SL_2(\Z)$ such that $W [a] = [L]$, if $W^{-1} H W$ may be expressed as a special form
\[ P^J A^{n(k)} \dots B^{n(2)} A^{n(1)}, \mbox{ where } J \in \{+1,-1\} \mbox{ and } |n(i)| \geq 2 \mbox{ for } i = 2, \dots, k\]
or
\[ P^J B^{n(k)} \dots B^{n(2)} A^{n(1)}, \mbox{ where } J \in \{0,2\} \mbox{ and } |n(i)| \geq 2 \mbox{ for } i = 2, \dots, k\]
then there is a properly embedded essential surface in $M-N(L)$ of type $C(J;n(k), \dots, n(1))$.
\end{prop}

\begin{proof}
Let $\omega$ be the homeomorphism of $T$ associated to $W$.  Thus $W[a] = [\omega(a)]$.  Then $\omega$ naturally extends to a bundle equivalence
\[\omega \colon M' = T \times I/_{\omega^{-1} \circ \eta \circ \omega} \to M = T \times I/_\eta \]
such that $\omega(a) = L$.
If $W^{-1} H W$ may be expressed as the special form $P^J \dots B^{n(2)} A^{n(1)}$ then 
\[ \omega^{-1} \circ \eta \circ \omega = \phi^J \dots \beta^{n(2)} \alpha^{n(1)}. \]
Thus 
\[ M' =  T \times I/_{\omega^{-1} \circ \eta \circ \omega} = T \times I/_{\phi^J \dots \beta^{n(2)} \alpha^{n(1)}}, \]
and so $M'$ contains the twisted surface $C(J; n(k), \dots, n(2), n(1))$.  Since $P^J \dots B^{n(2)} A^{n(1)}$ is a special form, $|n(i)| \geq 2$ for $i = 2, \dots, k$.  Hence by Proposition~\ref{prop:2.5.1} the surface $C(J; n(k), \dots, n(2), n(1))$ is essential in the complement of $a \times \{0\}$.  Therefore $\omega(C(J; n(k), \dots, n(2), n(1)))$ is an essential surface in $M-N(L) = \omega(M') - N(\omega(a))$ of type $C(J; n(k), \dots, n(2), n(1))$.
\end{proof}

Together, Propositions~\ref{prop:uniqueness} and~\ref{prop:existence} imply
\begin{thm}\label{thm:sfceandspecform}
Let $M$, $L$, and $H$ be as above.
The essential surfaces in $M - N(L)$ with the type of a twisted surface are in one-to-one correspondence with special forms of conjugates $W^{-1} H W$ of $H$ where $W \in \SL_2(\Z)$ such that $W [a] = [L]$. 
\end{thm}

\section{Structure of surfaces}\label{sec:surfacestructure}
 
This section closely follows \cite[\S3]{cjr:isioptb}.  Many of their arguments are easily modified to accommodate the presence of an essential level knot $L$ in a once-punctured torus bundle $M$.  Parts that do not involve $L$ will be simply cited. 

Let $S$ be a properly embedded surface in $M$ that is disjoint from $L$ and essential in $M-N(L)$.  We will say that $S$ is in {\em general position} provided that 
\begin{enumerate}
\item each component of $\bdry S$ is either contained in a fiber or is transverse to every fiber,
\item the projection $p:M \to S^1$ restricts to a Morse function on the interior of $S$ having distinct critical values different from $p(L)$,
\item among all surfaces isotopic to $S$ in $M-N(L)$ and satisfying $(1)$ and $(2)$, $S$ has the minimal number of index $0$ or $2$ critical points.
\end{enumerate}
We may assume that $S$ has been moved by an isotopy with support outside of a neighborhood of $L$ to be in general position.  The level sets of $p|_S$ are the intersections of $S$ with the fibers of $M$.

Let $x$ be a critical point.  We use the following terms as defined in \cite[\S3]{cjr:isioptb}: {\em level arcs}, {\em level sets}, {\em critical neighborhood} of $x$, and {\em upper} and {\em lower level sets of $x$}.  We then have the following:
\begin{lemma}\label{lem:lemmasof3}
$\mbox{}$
\begin{itemize}
\item[(a)]{\em(\cite[Lemma 3.1.1]{cjr:isioptb})}
Each level arc of $S$ is essential in the fiber containing it.
\item[(b)]{\em(\cite[Lemma 3.2.1]{cjr:isioptb})}
Either $S$ meets every non-critical fiber only in arcs or $S$ meets every non-critical fiber only in simple closed curves.
\item[(c)]{\em(\cite[Lemma 3.2.2]{cjr:isioptb})}
If both the upper and lower level sets of a critical neighborhood of a critical point contain an essential closed curve, then these curves are isotopic.
\item[(d)]{\em(\cite[Lemma 3.2.3]{cjr:isioptb})}
If the lower level set contains two arcs then
\begin{itemize}
\item[(i)] they are parallel, and
\item[(ii)] the upper level set is obtained by a band sum across the annulus component of the complement.
\end{itemize}
\end{itemize}
\end{lemma}

\begin{proof}
All four lemmas follow almost exactly as in \cite{cjr:isioptb}.

One needs for the proof of (b) that $M-N(L)$ is irreducible.  This follows from the irreducibility of $M$ and that $L$ is non-trivial in its fiber.

Since critical points of $S$ occur away from $L$, arguments involving critical neighborhoods are unchanged.
\end{proof}

For the proof of \cite[Theorem 3.3.1]{cjr:isioptb}, the authors employ a lemma due to Haken \cite{haken}.  For the upcoming proof of Theorem~\ref{thm:surface_type}, we must alter the lemma to accommodate the presence of $L$.

\begin{lemma}\label{lem:TxI}
Let $T$ be a punctured torus and $L$ be an essential simple closed curve in the fiber $T \times \{\frac{1}{2}\}$ of $T \times I$.  Let $R$ be a properly embedded, connected, incompressible surface in $T \times I - N(L)$ such that each component of $\bdry R$ is contained in either $\bdry T \times I$, $T \times \{0\}$, or $T \times \{1\}$ and, in the latter two cases, parallel in its fiber to $\bdry T \times I$.  Then $R$ is either an annulus, a once-punctured torus, or a torus parallel to $\bdry N(L)$.
\end{lemma}

\begin{proof}

$R$ is disjoint from $L$ in $T \times I$.  If $R$ is incompressible in $T \times I$, then it follows from \cite{haken} that $R$ is either an annulus or a once-punctured torus.  Therefore assume $R$ is compressible in $T \times I$ yet incompressible in the complement of $L$.  Thus any compressing disk must then intersect $L$.  

We may assume that $L = a \times \{\frac{1}{2}\}$.  Then let $A_a$ be the annulus $a \times I$, and let $A_b$ be the annulus $b \times I$.  Note that $L$ is the core of $A_a$ and intersects $A_b$ once.  

By standard arguments, we may assume that $R$ has been isotoped rel-$\bdry$ to minimize both $|A_a \cap R|$ and $|A_b \cap R|$.  It follows that $A_a \cap R$ is a collection of simple closed curves on $A_a$ that are parallel to the core of $A_a$ (and hence $L$), and that $A_b \cap R$ is a collection of simple closed curves on $A_b$ that either bound disks in $A_b$ that intersect $L$ or are parallel to the core of $A_b$.

{\em Case 1.}  Assume there exists a curve in $A_b \cap R$ that bounds a disk in $A_b$.  Among the curves of $A_b \cap R$ that bound disks in $A_b$, let $c$ be the innermost.  Let $D_c$ be the disk bounded by $c$.  By minimality assumptions, we may assume $|c \cap A_a| = 2$. Thus there exists a curve $c' \in A_a \cap R$, such that $c \cap c' \neq \emptyset$. Note that $c'$ is innermost among $A_a \cap R$ in the sense that if $A_{c'}$ is the subannulus of $A_a$ between $L$ and $c'$, then $A_{c'} \cap R = c'$.  Consider the loop $\delta = \bdry N(c \cup c') \cap R$ and the subdisk $\Delta$ of $\bdry N(D_c \cup A_{c'})$ that it bounds.  Since $\Delta$ is disjoint from $L$, $\delta$ must bound a disk, say $\Delta'$, in $R$.  Furthermore, note that $\Delta' \cap N(c \cup c') = \delta$.  Thus $R = \Delta' \cup N(c \cup c')$ which is a torus parallel to $\bdry N(L)$.

{\em Case 2.}  Assume that there are no curves of $A_b \cap R$ that bound disks in $A_b$.  Thus all curves of $A_b \cap R$ are parallel to the core of $A_b$.  

Among the compressing disks for $R$ in $T \times I$, let $D$ be one such that $|D \cap A_a|$ is minimal.  Thus $D \cap A_a$ contains no simple closed curve, only arcs.  Note that $D \cap A_a \neq \emptyset$ since $D \cap L \neq \emptyset$.  Due to the minimality assumptions on $|D \cap A_a|$ and $|R \cap A_a|$, every arc of $D \cap A_a$ interesects $L$ transversely. 

 Let $d$ be an outermost arc on $D$ of $D \cap A_a$, and let $D'$ be the (outermost) subdisk of $D$ cut off by $d$.  Let $d'$ be the arc of $\bdry D$ so that $\bdry D' = d' \cup d$.  The endpoints of $d$ lie on curves of $A_a \cap R$, say $c$ and $c'$, that are each adjacent to $L$.  Let $A_{c c'}$ be the subannulus in $A_a$ between $c$ and $c'$.  Its interior is disjoint from $R$.
 
Consider $\delta = \bdry (N(c \cup c' \cup d') \cap R)$ and the subdisk $\Delta$ of $\bdry N(A_{c c'} \cup D')$ that it bounds.  Since $\Delta \cap L = \emptyset$, $\delta$ must bound a disk, say $\Delta'$, in $R$.  Therefore we have the annulus $\Delta' \cup (N(c \cup c' \cup d') \cap R)$ which intersects $A_b$ in an arc $r$ of some component $c''$ of $A_b \cap R$.  Since $r$ connects $c \cap A_b$ to $c' \cap A_b$, $c''$ intersects $A_a \cap A_b$ at least twice.  However, since the curves of $A_b \cap R$ are parallel to the core, by the minimality assumptions on $|A_a \cap R|$, $|c'' \cap (A_a \cap A_b)| = 1$.  Thus we have a contradiction.
\end{proof}

\begin{thm} \label{thm:surface_type} Let $M$ be a once-punctured torus bundle such that the characteristic class of $M$ does not have trace $2$.  Let $L$ be an essential level curve in $M$.  If $(S, \bdry S) \subset (M-N(L), \bdry M)$ is an essential connected surface in $M-N(L)$ such that $\bdry S$ is transverse to the induced fibration on $\bdry M$, then $S$ is 
a surface of type $C(J;n(k), \dots, n(1))$.
\end{thm}

\begin{remark}
One may prove a more general version of Theorem~\ref{thm:surface_type} analogous to \cite[Theorem~3.3.1]{cjr:isioptb} for surfaces in the complement of an essential level knot such as $L$.  For our purposes, however, we restrict attention to once-punctured torus bundles and surfaces as in the theorem.
\end{remark}

\begin{proof}
This proof largely follows the proof of \cite[Theorem~3.3.1]{cjr:isioptb}.  We sketch only the parts of their proof that need alterations for the proof at hand.

We assume $S$ has been isotoped into general position.  By Lemma~\ref{lem:lemmasof3}(b) (\cite[Lemma 3.2.1]{cjr:isioptb}) either $S$ meets every non-critical fiber only in arcs or $S$ meets every non-critical fiber in simple closed curves.

{\em Case 1.} $S$ meets every non-critical fiber only in simple closed curves.

First assume there exists a non-critical fiber $F$ not containing $L$ with $S$ intersecting $F$ in only inessential curves, i.e.\ trivial or \bdry-parallel curves.  By a further isotopy of $S$, we may assume that the number of such curves is minimal.  Thus either $S \cap F = \emptyset$ or every component of $S \cap F$ is parallel into $\bdry F$.

Splitting $M$ along $F$ yields the product $T \times I$ and induces a splitting of $S$.  Let $S'$ be the surface resulting from splitting $S$ along $F$.  $S'$ is incompressible in $T \times I - N(L)$, and it follows that every component of $S'$ is just like $R$ in Lemma~\ref{lem:TxI}.  Thus every component of $S'$ is either an annulus, a once-punctured torus, or a torus parallel to $\bdry N(L)$.
A component of $S'$ cannot be a torus parallel to $\bdry N(L)$ since then $S = S'$.  This contradicts the assumption that $S$ is essential.

By considering the placement of the boundary of a component of $S'$, the minimality condition imposed on $S \cap F$ forces two things.  Either $|S \cap F|=1$ or $S \cap F = \emptyset$.  If $|S \cap F| = 1$, then $S'$ is an annulus parallel to $\bdry T \times I$ and $S$ is a torus parallel to $\bdry M$.  If $S \cap F = \emptyset$, then $S$ is a once punctured torus isotopic to a fiber.  Neither of these may occur as $S$ is not essential in the former and $\bdry S$ is not transverse to the fibration in the latter.

Now assume every non-critical fiber $F$ contains an essential simple closed curve component of $S \cap F$.  From here, we follow \cite{cjr:isioptb} exactly to the conclusion that the characteristic class of $M$ fixes the isotopy class of an essential curve.  Hence the characteristic class of $M$ has trace $2$, contrary to  assumption.

{\em Case 2.}  $S$ meets every non-critical fiber in essential arcs.  

The fiber containing $L$ must have exactly one non-empty family of parallel level arcs.  This allows us to skip much of \cite[Theorem 3.3.1, Case 2]{cjr:isioptb}.

Let $F$ be a fiber in a neighborhood of the fiber containing $L$.  $S$ may only meet $F$ in one family of parallel arcs.  If $S$ has no critical points, then $S$ meets every fiber in just one family of parallel arcs.  The characteristic class of $M$ must therefore fix the isotopy class of an essential simple closed curve and hence have trace $2$.  By assumption, this is not the case.  Thus $S$ must contain a critical point.  

The proof completes just as the last paragraph of \cite[Theorem~3.3.1]{cjr:isioptb}.  It follows that $S$ is a surface of type $C(J;n(k), \dots, n(1))$. 
\end{proof}

\section{Framing}\label{sec:framing}

\subsection{Review}
Let us recall the discussion about framing in \cite[\S 6.2]{cjr:isioptb}.  We will use the definitions and notation (and its abuse) established there.  The following is paraphrased from that section:

Fix a base point $x \in \bdry T$ and let $a$ and $b$ be elements of $\pi_1(T,x)$ analogous to the $a$ and $b$ in Figure~\ref{fig:o-pt}.  
  Let $\Stab([a,b])$ be the subgroup of the automorphisms of $\pi_1(T,x)$ that stabilizes $[a,b]$ and therefore act as identity on $\bdry T$ (i.e. up to isotopy fixing $x$).   Elements $\gamma \in \Stab([a,b])$ correspond uniquely (i) to homeomorphisms $g$ of $(T,x)$ such that $g_* = \gamma$ and (ii) to oriented simple closed curves $t_\gamma$ swept out from $x$ on the boundary of $M = T \times I/_g$ such that $t_\gamma a t_\gamma^{-1} = \gamma(a)$ and $t_\gamma b t_\gamma^{-1} = \gamma(b)$.  Note $t_\gamma$ is transverse to the fibration, intersecting each fiber once.  Together with the oriented boundary of a fiber (which is analogous to the standard longitude of a knot in $S^3$), the curve $t_\gamma$ defines a basis, or framing, for $H_1(\bdry M)$.  Therefore we have the following definition. 

An element $\gamma$ of $\Stab([a,b])$ is a {\em framing} for a once-punctured torus bundle $M = T \times I/_g$ if $g_* = \gamma$.

If $g'$ is also a homeomorphism of $(T,x)$ that is isotopic (not necessarily fixing $x$) to $g$, then there exists a bundle equivalence $h \colon T \times I/_{g'} \to  T \times I/_g$.
As elements of $\pi_1(T,x)$, the corresponding framings $\gamma = g_*$ and $\gamma' = g'_*$ differ by a conjugation by $[a,b]^j$ for some $j \in \Z$, i.e. $\gamma^{-1} h_* \gamma' h_*^{-1}$ is a conjugation by $[a,b]^j$.  This $j$ describes the difference in the number of times $t_\gamma$ and $t_{\gamma'}$ spin around the boundary of the fiber.
This number $j \in \Z$ is called the {\em transition index}  between $\gamma$ and $\gamma'$.

Given a framing for a once-punctured torus bundle, we may describe the boundary of an essential surface of type $C(J; n(k), \dots, n(1))$ in terms of this framing.  More precisely, we may associate a framing to the special form corresponding to the essential twisted surface and thereby understand the boundary curves of the twisted surface in terms of this framing.

The {\em standard framings} defined in \cite{cjr:isioptb} are as follows:
\[\begin{array}{ll}
\alpha\colon \begin{cases} a \to a \\ b \to b a^{-1} \end{cases}
&\mbox{ is standard for } A = \mat{1}{0}{-1}{1}, \\
\beta\colon \begin{cases} a \to a b \\ b \to b \end{cases}
&\mbox{ is standard for } B = \mat{1}{1}{0}{1}, \\
\phi\colon \begin{cases} a \to a b a^{-1} \\ b \to a^{-1} \end{cases} 
&\mbox{ is standard for } P = \mat{0}{1}{-1}{0}, \\
\psi\colon \begin{cases} a \to a b^{-1} \\ b \to b a b^{-1} \end{cases}
&\mbox{ is standard for } Q = \mat{1}{-1}{1}{0}.
\end{array}\]
Note that $\alpha$, $\beta$, $\phi$, and $\psi$ are being used duplicitously for elements of $\Stab([a,b])$ and homeomorphisms of $T$.  



Observe that $\Stab([a,b])$ is generated by the standard framings $\alpha$ and $\beta$ together with the framing
\[\delta\colon x \to [a,b] x [a,b]^{-1} \]
induced by a Dehn twist along a simple curve parallel to $\bdry T$.
Since $\alpha = \psi\phi$, $\beta = \phi\psi$, and $\delta = \phi^4$, $\Stab([a,b])$ is generated by the standard framings $\phi$ and $\psi$.  Indeed, we have the presentation
\[ \Stab([a,b])= \langle \phi, \psi | \phi^2 \psi^3 = 1 \rangle. \]
Note that $\phi^2$ and $\psi^3$ are both in the center of $\Stab([a,b])$.  
Moreover, if $j$ is the transition index between two framings $\xi$ and $\zeta$, then 
\[ \xi^{-1} \mu \zeta \mu^{-1} = \phi^{4j} \]
for some framing $\mu$.

\subsection{Exponent sums}

Let $\omega$ be a word on the set of letters $\Lambda$, say $w = \prod_{i=1}^k \lambda_i^{e_i}$ where $\lambda_i \in \Lambda$ and $e_i \in \Z$.  Define the {\em exponent sum} of the word $\omega$ to be $E(\omega) = \sum_{i=1}^k e_i$.  Define the {\em exponent sum of the letter $\lambda$} for the word $\omega$ to be $E_\lambda(\omega) = \sum_{\lambda_i = \lambda} e_i$.

\begin{lemma}\label{lem:spelling}
Let $\omega$ be a word for an element of $\Stab[a,b]$ written in the letters $\phi$ and $\psi$.  Then $(E_\phi(\omega), E_\psi(\omega)) \mod (2,3)$ is independent of the spelling of $\omega$.
\end{lemma}

\begin{proof}
Let $\omega$ and $\omega'$ be two words in $\phi$ and $\psi$ representing the same element of $\Stab([a,b])$.  The passage between the two words occurs by the following moves:
\begin{itemize}
\item[$(i)$] insertion or deletion of an adjacent pair of a letter and its inverse and
\item[$(ii)$] application of the relation $\phi^2 \psi^3 = 1$.
\end{itemize}
Assume $\omega$ and $\omega'$ differ by just one of these moves.  Since $(i)$ only introduces or removes canceling pairs, it does not change the exponent sums at all.  For $(ii)$, assume $\omega'$ is obtained from $\omega$ by replacing $\phi^2 \psi^3$ with $1$.  Then 
\[(E_\phi(\omega'), E_\psi(\omega')) = (E_\phi(\omega)-2, E_\psi(\omega)-3) = (E_\phi(\omega),E_\psi(\omega))-(2,3)\]
Thus 
\[(E_\phi(\omega'), E_\psi(\omega')) = (E_\phi(\omega), E_\psi(\omega)) \mod (2,3). \] 
\end{proof}

\begin{lemma} \label{lem:framing}
Let $\zeta$ be a framing for $Z \in \SL_2(\Z)$, and let $\xi$ be a framing for $X\in \SL_2(\Z))$ conjugate to $Z$.  Let $\omega$ be a word for $\zeta^{-1} \xi$ written in the letters $\phi$ and $\psi$.  Then 
\[ (E_\phi(\omega), E_\psi(\omega)) = (4j, 0) \mod (2,3)\]
where $j$ is the transition index between $\zeta$ and $\xi$.
\end{lemma}

\begin{proof}
Let $j$ be the transition index between $\zeta$ and $\xi$.  Thus $\zeta^{-1} \mu \xi \mu^{-1} = \phi^{4j}$ for some framing $\mu$ associated to $U \in \SL_2(\Z)$ such that $Z^{-1} U X U^{-1} = I$.  Let $\omega_\zeta$, $\omega_\xi$, and $\omega_\mu$ be words in $\phi$ and $\psi$ for $\zeta$, $\xi$, and $\mu$ respectively.  Thus we may write
\[ \omega_\zeta^{-1} \omega_\mu \omega_\xi \omega_\mu^{-1} = \phi^{4j}. \]

Since $\zeta^{-1} \mu \xi \mu^{-1} = \phi^{4j}$, the words $\omega_\zeta^{-1} \omega_\mu \omega_\xi \omega_\mu^{-1}$ and $\phi^{4j}$ represent the same element of $\Stab([a,b])$.  By Lemma~\ref{lem:spelling}
\[
\begin{array}{rll}
(E_\phi(\omega_\zeta^{-1} \omega_\mu \omega_\xi \omega_\mu^{-1}),E_\psi(\omega_\zeta^{-1} \omega_\mu \omega_\xi \omega_\mu^{-1})) &= (E_\phi(\phi^{4j}), E_\psi(\phi^{4j})) &\mod (2,3) \\
&= (4j,0) &\mod (2,3).
\end{array}
\]
Moreover, since $E_\phi(\omega_\mu) = -E_\phi(\omega_\mu^{-1})$ and $E_\psi(\omega_\mu) = -E_\psi(\omega_\mu^{-1})$, we have
\[ (E_\phi(\omega_\zeta^{-1} \omega_\mu \omega_\xi \omega_\mu^{-1}), E_\psi(\omega_\zeta^{-1} \omega_\mu \omega_\xi \omega_\mu^{-1})) = (E_\phi(\omega_\zeta^{-1} \omega_\xi), E_\psi(\omega_\zeta^{-1} \omega_\xi)).\]
Because $\omega_\zeta^{-1} \omega_\xi$ is a word for $\zeta^{-1} \xi$, the conclusion of the lemma follows.
\end{proof}

Assume $P^J \dots B^{n(2)} A^{n(1)}$ is a special form for an element $X$ of $\SL_2(\Z)$.  This corresponds to an essential twisted surface in the complement of an essential level knot in a certain once-punctured torus bundle $M$ with characteristic class $[X]$.  Consider the standard framing $\xi = \phi^J \dots \beta^{n(2)} \alpha^{n(1)}$ corresponding to the specific form for $X$.  Assume a framing $\zeta$ has been chosen for $M$, associated to $Z \in \SL_2(\Z)$.  Thus $[Z] = [X]$ and  $Z^{-1} U X U^{-1} = I$ for some $U \in \SL_2(\Z)$ whereas $\zeta^{-1} \mu \xi \mu^{-1} = \phi^{4j}$ where $\mu$ is any framing for $U$ and $j \in \Z$ is the transition index.

\begin{prop}\label{prop:meridian_conditions}
Continuing with the above notation, let $t_\zeta$ be the curve on $\bdry M$ associated to $\zeta$.
Assume $Z$ may be written as a word $W_Z$ in $A$ and $B$ so that the corresponding standard framing is $\zeta$.
The boundary components of $C(J; n(k), \dots, n(2), n(1))$ are isotopic to $t_\zeta$ if and only if $J=0$ and $E(W_Z) = \sum_{i=1}^{k} n(i)$.
\end{prop}

\begin{remark}
Here one might care to consider $M$ as the exterior of some genus one fibered knot with $t_\zeta$ as the meridian and the boundary of the fiber as the standard longitude.
\end{remark}

\begin{proof}
A boundary component of $C(J; n(k), \dots, n(2), n(1))$ intersects a fiber of $M$ only once if and only if $J = 0$.  The transition index between $\xi$ (the framing corresponding to the twisted surface) and $\zeta$ measures how many times the boundary component wraps longitudinally.  Thus the boundary components of $C(J; n(k), \dots, n(2), n(1))$ are isotopic to $t_\zeta$ if and only if $J = 0$ and the transition index between $\xi$ and $\zeta$ is zero.

Recall that the standard framings for $A$ and $B$ are $\alpha = \psi \phi$ and $\beta = \phi \psi$ respectively.  We may therefore consider words in $\alpha$ and $\beta$ as actually words in $\phi$ and $\psi$.  It then follows that  $E_\phi$ and $E_\psi$ are equal on any word written in $\alpha$ and $\beta$. 
Indeed, given a word $\omega$ written in $\alpha$ and $\beta$, 
\[ (E_\phi(\omega), E_\psi(\omega)) = E(\omega)\cdot(1,1) \] 
where $E$ is the exponent sum over the letters $\alpha$ and $\beta$.

Let $\omega_\zeta$ be the standard framing for $W_Z$ written as a word in $\alpha$ and $\beta$.  Then since $\xi = \phi^J \dots \beta^{n(2)} \alpha^{n(1)}$ and $J=0$, let $\omega_\xi$ be the word $\beta^{n(k)} \dots \beta^{n(2)} \alpha^{n(1)}$.

Thus we have $\zeta^{-1} \xi$ written as the word $\omega_\zeta^{-1} \omega_\xi$ in $\alpha$ and $\beta$.  Therefore
\begin{multline*}(E_\phi(\omega_\zeta^{-1} \omega_\xi), E_\psi(\omega_\zeta^{-1} \omega_\xi)) = E(\omega_\zeta^{-1} \omega_\xi)\cdot(1,1) \\ 
= (E(\omega_\zeta^{-1} \omega_\xi),E(\omega_\zeta^{-1} \omega_\xi)) \mod (2,3). \end{multline*}
Since we need the transition index to be zero, we require 
\[ (E(\omega_\zeta^{-1} \omega_\xi), E(\omega_\zeta^{-1} \omega_\xi)) = (0,0) \mod (2,3).\]
In other words, it must be that
\[ (E(\omega_\zeta^{-1} \omega_\xi), E(\omega_\zeta^{-1} \omega_\xi)) - N\cdot(2,3) = (0,0)\]
for some $N \in \Z$.  This can only happen in the case that $E(\omega_\zeta^{-1}\omega_\xi) = 0$ and $N=0$.  Therefore $E(\omega_\zeta) = E(\omega_\xi)$.

Since $E(\omega_\zeta) = E(W_Z)$ and $E(\omega_\xi) = \sum_{i=1}^k n(i)$, the conclusion of the proposition follows.
\end{proof}

\section{Algorithms}


\begin{algor}\label{algor:main_algorithm}
Let $M$ be a once-punctured torus bundle given by an element $H \in \SL_2(\Z)$.  Let $L$ be an essential simple closed curve on the fiber $T \times \{0\}$ of $M$ given by the ordered pair $(x, y) \in \Z^2$.  The following steps give a procedure to list all essential surfaces of type $C(J; n(k), \dots, n(1))$ in $M - N(L)$.
\begin{enumerate}
\item[{\em Step 1.}]  Choose a change of basis matrix $W \in \SL_2(\Z)$ so that $W \begin{pmatrix} 1 \\ 0 \end{pmatrix} = \begin{pmatrix} x \\ y \end{pmatrix}$.  

\item[{\em Step 2.}]  Let $X(N) = (W A^N)^{-1} H (W A^N)$.

\item[{\em Step 3.}]  List all $N$ such that $X(N) = \begin{pmatrix} p(N) & r(N) \\ q(N) & s(N) \end{pmatrix}$ has 
\[ |p(N)| < |q(N)| \]

\item[{\em Step 4.}]  For each $N$ listed in Step 3, 
\begin{itemize}
\item obtain the minimal continued fraction expansions $[a_1, a_2, \dots, a_{k-1}]$ for $p(N)/q(N)$  
\item find $a_k$ such that $[a_1, a_2, \dots, a_{k-1}, a_k]$ is a continued fraction expansion for $r(N)/s(N)$.
\end{itemize}
If $k$ is odd, then for either $J = +1$ or $J = -1$,
\[ P^J A^{a_1} B^{a_2} \dots A^{a_{k-2}} B^{a_{k-1}} A^{a_k} = X(N). \]
If $k$ is even, then for either $J = 0$ or $J = 2$, 
\[ P^J B^{a_1} A^{a_2} \dots A^{a_{k-2}} B^{a_{k-1}} A^{a_k} = X(N). \] 

\item[{\em Step 5.}]  For each $N$ in Step 3, list the expressions for $X(N)$ obtained in Step 4 rewriting $a_i$ as $n(k-i+1)$.

\item[{\em Step 6.}]  For every expression $P^J  \dots B^{n(2)} A^{n(1)}$ listed in Step 5, there is an essential surface of type $C(J; n(k), \dots, n(1))$ in $M - N(L)$.
\end{enumerate}
\end{algor}

\begin{proof}
Given a change of basis matrix $W$ such that $W[a] = [L]$, all other such matrices are of the form $W A^N$ for some integer $N$.  Proposition~\ref{prop:existence} then implies that each expression of $X(N) = (W A^N)^{-1} H (W A^N)$ as a special form corresponds to an essential surface in $M - N(L)$ with the type of an essential twisted surface.  Steps 1 and 2 set up $X(N)$, and the remaining steps compute the special forms for all the $X(N)$.

By Lemma~\ref{lem:specialformMCFE} if $X(N) = \mat{p(N)}{q(N)}{r(N)}{s(N)}$ has a special form, then $p(N)/q(N)$ has a MCFE.  This then implies that $|p(N)/q(N)| < 1$ and hence $|p(N)| < |q(N)|$.   Since $X(N) = A^{-N} (W^{-1}HW) A^N$, $q(N)$ does not depend on $N$.  Thus the list of Step 3 is a finite list.

For each $N$ such that $p(N)/q(N)$ has a MCFE $[a_1, a_2, \dots, a_{k-1}]$ we may solve for $a_k$ as done in Lemma~\ref{lem:matrix_cfe2} so that $[a_1, \dots, a_{k-1}, a_k] = r(N)/s(N)$.  Lemma~\ref{lem:specialformMCFE} gives the associated special form for $X(N)$.   Step 4 lists these special forms for $X(N)$.

Step 5 collects all the special forms for $X(N)$ for all $N$.  Step 6 then follows from Proposition~\ref{prop:existence}.
\end{proof} 

 Given a framing $\zeta$ for $M$ we extend Algorithm~\ref{algor:main_algorithm} to describe the boundary curves of the essential surfaces of type $C(J; n(k), \dots, n(1))$ in $M - N(L)$ in terms of the given framing.  
This may be done by a direct computation which has been previously described in \cite[\S6.2]{cjr:isioptb}.  Lemma~\ref{lem:framing}, however, allows us to simplify this computation.

\begin{algor}[Algorithm~\ref{algor:main_algorithm} Continued: Framing]\label{algor:framing_algorithm}
Let $\zeta$ be any given framing for $M$.
For each special form $P^J \dots B^{n(2)} A^{n(1)}$ listed in Step 5 of Algorithm~\ref{algor:main_algorithm}, compute the coordinates of the boundary of the associated essential surface as follows:

\begin{itemize}
\item[{\em Step 7.}]  Take $\xi$ to be the standard framing associated to the special form.

\item[{\em Step 8.}]  Write $\zeta^{-1}\xi$ as a word $\omega$ in the letters $\phi$ and $\psi$.

\item[{\em Step 9.}]  Find $j \in \Z$, the transition index,  such that
\[ (E_\phi(\omega), E_\psi(\omega)) = (4j, 0) \mod (2,3). \]
\end{itemize}

With respect to the framing $\zeta$, the boundary components of the essential surface of type $C(J; n(k), \dots, n(1))$ have coordinates $\langle 1, -j \rangle$, $\langle 2, 1-j \rangle$, $\langle 4, 1-j \rangle$, or $\langle 4, -1-j \rangle$ if $J = 0$, $2$, $+1$, or $-1$ respectively.
 \end{algor}
 
 Here the coordinates $\langle p,q \rangle$ of a curve on $\bdry M$ mean that, up to sign, the curve is homologous to $p[t_\zeta] + q [\bdry T \times \{0\}] \in H_1(\bdry M)$.

\begin{proof}
In Step 8 we may write $\zeta^{-1}\xi$ in terms of $\phi$ and $\psi$ by using the normal form for the structure of $\Stab([a,b])$ as a free product with amalgamation.   In Step 9 we find the transition index via Lemma~\ref{lem:framing}.  \cite[Table~1]{cjr:isioptb} lists the boundary curves for twisted surfaces in terms of their corresponding standard framings.  The transition index tells us how many times around the boundary of the fiber (the longitude) the boundary of the twisted surface wraps with respect to the given framing $\zeta$.
\end{proof}

\section{Passing between closed surfaces and surfaces with boundary}\label{sec:closedsurfaces}

Let $\widehat{M}$ be a $3$--manifold containing knots $L$ and $K$.
Let $\widehat{S}$ be an essential closed (orientable, connected) surface $\not \cong S^2$ in $\widehat{M}-N(L)$.  Assume $\widehat{S}$ has been chosen among surfaces in its isotopy class in $\widehat{M}-N(L)$ so that $|K \cap \widehat{S}|$ is minimized.  Let $S = \widehat{S} - N(K) \subset \widehat{M}-N(K \cup L)$.  Note that $\bdry S$ is a collection of meridional curves on $\bdry N(K)$.

\begin{lemma} \label{lem:ess_with_bdry}
$S$ is essential in $\widehat{M}-N(K \cup L)$.
\end{lemma}

\begin{proof}
We must show that $S$ is incompressible, $\bdry$-incompressible, and not $\bdry$-parallel.  

Assume $S$ is compressible in $\widehat{M}-N(K \cup L)$.  Let $D$ be a compressing disk for $S$.  Since $\widehat{S}$ is incompressible, $\bdry D$ must bound a disk $E \subset \widehat{S} \subset \widehat{M}-N(L)$.  $E$ must intersect $K$ since otherwise $\bdry D$ would bound the disk $E$ in $S$ contradicting that $D$ is a compressing disk.  Let $\widehat{S}'= (\widehat{S} - N(E)) \cup D$.  Since $\widehat{M}-N(L)$ is irreducible, $\widehat{S}'$ is isotopic to $\widehat{S}$.  But then $|\widehat{S} \cap K| > |\widehat{S}' \cap K|$ contradicting the minimality assumption.

Assume $S$ is $\bdry$-parallel.   Then $S$ is either a torus parallel to $\bdry N(L)$, a torus parallel to $\bdry N(K)$, or an annulus parallel into $\bdry N(K)$.  If $S$ is parallel to $\bdry N(L)$ then $S = \widehat{S}$, and so $\widehat{S}$ is $\bdry$-parallel in $\widehat{M}-N(L)$.  If $S$ is parallel to $\bdry N(K)$, then $S = \widehat{S}$, and so $S$ is compressible in $\widehat{M}-N(L)$.  If $S$ is an annulus parallel into $\bdry N(K)$, then since $\bdry S$ is a collection of meridional curves, $S \cong S^2$ in $\widehat{M}-N(L)$.  These all contradict our assumptions on $S$.

Assume $S$ is $\bdry$-compressible, incompressible, and not $\bdry$-parallel.  Since the boundary components of $\widehat{M}-N(K \cup L)$ are tori, $S$ must be an annulus.  Since $\bdry S$ is a collection of meridional curves, $S \cong S^2$ in $\widehat{M}-N(L)$ contrary to our assumption.
\end{proof}

In general the converse of Lemma~\ref{lem:ess_with_bdry} does not hold true.  Nevertheless, the characterization in Proposition~\ref{prop:2.5.1} generalizes for capped off twisted surfaces.

 Let $J, n(1), n(2), \dots, n(k)$, $L$ and $M$ be as usual.    Let $\widehat{M}$ be the manifold constructed by attaching a solid torus to $M$ so that the boundary curves of $C(J; n(k), \dots, n(2), n(1))$ are identified to meridional curves in the solid torus.  Let $\widehat{C}(J; n(k), \dots, n(2), n(1))$ be the closed surface in $\widehat{M}$ obtained by attaching disks to the boundary curves of $C(J; n(k), \dots, n(2), n(1))$.

\begin{thm}
[cf.\ Remark 2.5.2, \cite{cjr:isioptb}]
\label{thm:2.5.2}
  The closed surface $\widehat{C}(J; n(k), \dots, n(2), n(1))$ is incompressible in $\widehat{M}-N(L)$ if and only if $|n(i)| \geq 2$ for $i = 2, \dots, k$.  
\end{thm}

\begin{proof}
Let $K$ be the core of the attached solid torus so that $\widehat{M} - N(K) = \widehat{M}_K = M$.  Then $\widehat{C}(J; n(k), \dots, n(2), n(1)) - N(K) = C(J; n(k), \dots, n(2), n(1))$.  Let $C = C(J; n(k), \dots, n(2), n(1))$ and $\widehat{C} = \widehat{C}(J; n(k), \dots, n(2), n(1))$.  Note that $C$ is separating in $M$ and $\widehat{C}$ is separating in $\widehat{M}$. 

If $\widehat{C}$ is incompressible in $\widehat{M}_L$, then by Lemma~\ref{lem:ess_with_bdry}, $C$ is essential in $\widehat{M}_{K \cup L} = M_L$.  By Proposition~\ref{prop:2.5.1}, $|n(i)| \geq 2$ for $i \geq 2$.

Now assume $|n(i)| \geq 2$ for $i \geq 2$.  By Proposition~\ref{prop:2.5.1}, $C$ is essential in $M_L$.  For contradiction, assume $\widehat{C}$ is not essential in $\widehat{M}_L$.  Therefore there exists a compressing disk for $\widehat{C}$.

With $C \subset M$, express $M$ as
\[  M = T \times I /_{\alpha^{n(1)}} T \times I /_{\beta^{n(2)}} \dots T \times I/_{\phi^J \circ \gamma^{n(k)}}\]
where $\gamma = \beta$ if $k$ and $J$ are even and $\gamma = \alpha$ if $k$ and $J$ are odd.  If $J = 0$ or $2$, then each block meets $C$ in a single disk as in \cite[Figure 4]{cjr:isioptb}.  If $J = \pm 1$, then each block meets $C$ in two disks which are a parallel, each of which individually appears as in \cite[Figure 4]{cjr:isioptb}.  Let $F$ be the union of the fibers along which the blocks of $M$ are joined.

Consider the family of all compressing disks for $\widehat{C}$ in 
$\widehat{M}-N(L)$ that are transverse to $F$.  If this family is non-empty then there exists a member $D$ for which $(D \cap K, D \cap F)$ is minimized lexicographically.  Since $C$ is essential in $M-N(L)$, $D \cap K \neq \emptyset$.

 Let $P = D - N(K) \subset M-N(L)$ be the punctured disk of which the one boundary component on $C$ is $\bdry D$ and all other boundary components are on $\bdry N(K)$ parallel to $\bdry C$.  $P$ is incompressible in $M-N(L)$ due to the minimality assumption on $D$.  Since $D \cap K \neq \emptyset$, $\bdry P \cap \bdry M \neq \emptyset$.

\begin{claim} 
There cannot be any simple closed curve components of $P \cap F$ that bound disks in $M-N(L)$.
\end{claim} 

\begin{proof}
Since $F$ is incompressible, a simple closed curve component of $P \cap F$ must bound a disk in $F$.  An innermost such curve bounds a disk $\Delta \subset F$ with interior disjoint from $P$.  Chopping $P$ along $\Delta$ produces two planar surfaces.  Let $P'$ be the one with the boundary component $\bdry D$.  In $\widehat{M}-N(L)$, $P'$ caps off to a compressing disk $D'$ for $\widehat{C}$ with $\bdry D'  = \bdry D$.  However $|D' \cap K| \leq |D \cap K|$ and $|D' \cap F| < |D \cap F|$ contradicting the minimality assumption on $D$.
\end{proof}

{\em Part I. $J = 0$ or $2$}

Let $X$ be a block of $M$ with $\bdry D \cap X \neq \emptyset$.  Cutting $X$ along $C$ yields two solid tori as shown in \cite[Figure 6]{cjr:isioptb}.  Consider one of the solid tori which has non-empty intersection with the interior of $N(\bdry D) \cap P$.   Observe that the boundary of this solid torus meets $C$ in a disk, $\bdry M$ in two disks which we together label $B$, and $F$ in a disk $F_d$ and an annulus $F_a$.  

Suppose that $\sigma$ is an arc component of $P \cap F_d$.  
\begin{claim}\label{claim:arccase1}
Both end points of $\sigma$ cannot be on the same component of $B \cap F_d$.  
\end{claim}
\begin{proof}
Assume otherwise.  Thus $\sigma$ is an arc with both end points in $\bdry P - \bdry D \subset \bdry M$.  Among all the components of $P \cap F_d$ that have both end points on the same component of $B \cap F_d$, assume $\sigma'$ is outermost on $F_d$.  Thus there is a disk $\Delta \subset F_d$ such that $\bdry \Delta$ consists of two arcs, $\sigma'$ and $\delta'$ where $\delta' \subset B \cap F_d$ and $\Delta \cap P = \sigma'$.  Hence $\delta'$ is an arc on $\bdry M$ connecting two distinct components of $\bdry P$ on $\bdry M$.  


In $\widehat{M}$, $\delta'$ is just an arc on $K$.  Thus $\Delta$ guides an isotopy of $(D, \bdry D) \subset (\widehat{M}-N(L), \widehat{C})$ through $K$, reducing $|D \cap K|$.  This contradicts the minimality assumption.  (Observe that if $J=2$ the endpoints of $\sigma$ lie on components of $\bdry P$ that each intersect $F_d$ a second time on the other component of $B \cap F_d$.  The isotopy guided by $\Delta$ causes the arcs emanating from these other endpoints to connect.)
\end{proof}

\begin{claim}\label{claim:arccase2}
Both end points of $\sigma$ cannot be on the same component of $C \cap F_d$.  
\end{claim}
\begin{proof}
Assume otherwise.  Thus $\sigma$ is an arc with both end points in $\bdry D \subset C$.  Among all the components of $P \cap F_d$ that have both end points on the same component of $C \cap F_d$, assume $\sigma'$ is outermost on $F_d$.  Thus there is a disk $\Delta \subset F_d$ such that $\bdry \Delta$ consists of two arcs, $\sigma'$ and $\delta'$ where $\delta' \subset C \cap F_d$ and $\Delta \cap P = \sigma'$.

In $\widehat{M} - N(L)$, one may chop $D$ along $\Delta$ to produce two disks with boundary on $\widehat{C}$.  Both of these disks have fewer intersections with $F$ and no more intersections with $K$ than $D$.  Since $\bdry D$ is essential on $\widehat{C}$, at least one of their boundaries must be essential in $\widehat{C}$.  But then this disk is a compressing disk for $\widehat{C}$ that contradicts the minimality assumption on $D$.
\end{proof}

\begin{claim}\label{claim:arccase3}
$\sigma$ cannot have one end point contained in $B$ and the other contained in $C$.  
\end{claim}
\begin{proof}
Assume otherwise.
On $P$, $\sigma$ must connect $\bdry D$ with another boundary component of $P$.  In light of Claims~\ref{claim:arccase1} and \ref{claim:arccase2}, we may assume $\sigma'$ is an arc of $P \cap F_d$ with one end point in $B$ and one end point in $C$ an arc that is outermost on $F_d$ among all arcs of $P \cap F_d$.  Thus there is a disk $\Delta \subset F_d$ such that $\bdry \Delta$ is the union of three consecutive arcs $\sigma'$, $r' \subset C$, and $\delta' \subset B$, and $\Delta \cap P = \sigma'$.  $\delta'$ connects a component of $\bdry P$ and $\bdry C$ on $\bdry M$.  

As in Claim~\ref{claim:arccase1}, $\delta'$ is just an arc on $K$ in $\widehat{M}-N(L)$.  Thus $\Delta$ guides an isotopy of $(D, \bdry D) \subset (\widehat{M}-N(L), \widehat{C})$ through $K$, reducing $|D \cap K|$.  This contradicts the minimality assumption.
\end{proof}

\begin{claim}\label{claim:arccase4} 
$\sigma$ cannot connect both components of $C \cap F_d$
\end{claim}
\begin{proof}
Assume otherwise.  Due to Claims~\ref{claim:arccase1}, \ref{claim:arccase2} and \ref{claim:arccase3}, $P \cap F_d$ consists solely of a parallel family of arcs connecting the two arcs of $C \cap F_d$.  Then $P \cap F_d$  is disjoint from $B$.  Since $C$ is separating, this implies that $\bdry P \cap \bdry M = \emptyset$.  Thus $P = D$ contrary to assumptions.
\end{proof}

The fifth and final possible arc type for $\sigma$ connects the two components of $B \cap F_d$.  We need not address this case for our purposes.

Let $X'$ denote the block in $M$ that meets the block $X$ along the component of $F$ containing $F_a$.  $X'$ splits into two solid tori as $X$ did, and the same arguments of Claims~\ref{claim:arccase1}, \ref{claim:arccase2}, \ref{claim:arccase3} and \ref{claim:arccase4} for $X$ apply to $X'$.  Joining $X$ and $X'$ along the annulus $F_a$ forms a solid torus $V$ which meets $C$ in an annulus that wraps around the solid torus $|n(i)|$ times for some $i$. 

As a consequence of Claims~\ref{claim:arccase2}, \ref{claim:arccase3} and \ref{claim:arccase4}, $F_d \cap \bdry D= \emptyset$.  Therefore $\bdry D$ is contained in the annulus $C \cap V$ and is thus isotopic to the core of this annulus.

Consider the component $Q$ of $P \cap V$ that contains $\bdry D$.  Since $P$ is incompressible in $M-N(L)$, then $Q$ must be incompressible in $V - N(L)$.  Since all components of $\bdry Q$ are parallel on $\bdry V$, $Q$ is either a meridional disk of $V$ or a  boundary parallel annulus.  If $Q$ is a meridional disk then $Q=P=D$ contrary to assumptions.  Thus $Q$ is a boundary parallel annulus. 

One boundary component of $Q$, $\bdry D$, is contained in $C \cap V$ and the other is composed of two arcs on $\bdry M \cap V$ and two arcs on $F \cap V$.
There is a boundary compressing disk $\Delta$ for $Q$ in $V$ whose boundary is composed of a transverse arc $\sigma$ in $Q$, an arc in $C$ and an arc in $\bdry M$.

This arc $\sigma$ is like the arc in Claim~\ref{claim:arccase3}.  Similarly $\Delta$ guides an isotopy of $(D, \bdry D) \subset (\widehat{M}-N(L), \widehat{C})$ through $K$ that reduces $|D \cap K|$.  This contradicts the minimality assumption.

{\em Part II. $J = \pm 1$}

Again, let $X$ be a block of $M_{K}$ with $\bdry D \cap X \neq \emptyset$.  Cutting $X$ along $C$ yields two solid tori as shown in \cite[Figure 6]{cjr:isioptb} as well as one product disk ($ \cong D^2 \times I$) that gives the parallelism between the two disk components of $X \cap C$.

The arguments of Part I apply to these solid torus components.  Therefore we only need consider the case that the product disk has non-empty intersection with $\bdry D$.  The boundary of the product disk meets $C$ in two disks ($\cong D^2 \times \bdry I$), $\bdry M$ in four disks (labelled $B$), and $F$ in four disks. 

Each disk component of the intersection of the product disk with $F$ is a rectangle with two edges on $C$ and two edges on $B$.   Let $F_d$ be one of these components which has non-empty intersection with $\bdry D$.  Claims~\ref{claim:arccase1}, \ref{claim:arccase2}, \ref{claim:arccase3} and \ref{claim:arccase4} all then apply to the arcs of $P \cap F_d$.  Therefore $F_d \cap \bdry D = \emptyset$.  Thus $\bdry D$ must be contained in this product disk.  But this cannot occur as then $\bdry D$ would be contained in a disk of $C$ and hence bound a disk in $C$. 
\end{proof}

\begin{remark}
Theorem~\ref{thm:2.5.2} does not address whether an essential genus one twisted surface $C \subset M-N(L)$  might cap off to a torus $\widehat{C} \subset \widehat{M}-N(L)$ that is parallel to $\bdry N(L)$.  
\end{remark}



\section{Complements of essential surfaces}\label{sec:surfacecomplements}
  
Let $K$ be a genus one fibered knot in a closed $3$--manifold $\widehat{M}$.   Let $M =\widehat{M}-N(K)$ be the corresponding once-punctured torus bundle, and let $L$ be an essential level curve.

 We say that a surface $\widehat{S}$ in $\widehat{M}-N(L)$ has the {\em same type} as $\widehat{C}(J; n(k), \dots, n(2), n(1))$ if it is isotopic to $\widehat{C}(J; n(k), \dots, n(2), n(1))$ in $\widehat{M}-N(L)$.

\begin{lemma}\label{lem:n1leq1}
A closed surface of type $\widehat{C}(0; n(k), \dots, n(1))$ is incompressible in $\widehat{M}-N(L)$ but compressible in $\widehat{M}$ if and only if $|n(i)| \geq 2$ for $i \geq 2$  and $|n(1)| \leq 1$.
\end{lemma}

\begin{proof}
A surface $\widehat{S}$ of type $\widehat{C}(0; n(k), \dots, n(1))$ in $\widehat{M}-N(L)$ corresponds to a surface $S$ of type $C(0; n(k), \dots, n(1))$ in $M-N(L)$.  
By Proposition~\ref{prop:2.5.1} and Theorem~\ref{thm:2.5.2}, $S$ is essential and $\widehat{S}$ is incompressible if and only if $|n(i)| \geq 2$ for $i \geq 2$. 

By \cite[Proposition 2.5.1 and Remark 2.5.2]{cjr:isioptb}, $S$ is essential and $\widehat{S}$ is incompressible if and only if also $|n(1)| \geq 2$.  

Thus $\widehat{S}$ is compressible in $\widehat{M}$ and not $\widehat{M}-N(L)$ if and only if $|n(1)| \leq 1$ instead.
\end{proof}


\begin{remark}
If $\widehat{M} = S^3$, then since every closed surface must be compressible, a surface of type $\widehat{C}(0;n(k), \dots, n(2), n(1))$ that is incompressible in the complement of a level knot on the fiber of a trefoil or figure eight knot must have $|n(1)| \leq 1$ by the above lemma.
\end{remark}


\begin{prop}\label{prop:differences_in_surfaces}
Let $K$, $L$, and $M$ be as above.  Let $\widehat{S}$ be a closed essential surface in $\widehat{M}-N(L)$.
If $\widehat{S}$ is of type $\widehat{C}(0; n(k), \dots, n(2), n(1))$, then there exists an embedded annulus from $S$ to a curve on $\bdry N(L)$ of slope $-1/n(1)$.   In particular if $n(1) = \pm 1$, this slope is longitudinal; and if $n(1)=0$, this slope is meridional.
\begin{itemize}
\item If $|n(1)| \geq 2$, then $\widehat{S}$ is incompressible in $\widehat{M}$.
\item If $n(1) = \pm 1$, then $\widehat{S}$ bounds a handlebody in $\widehat{M}$.
\item If $n(1)=0$, then $\widehat{S}$ does not bound a handlebody in $\widehat{M}$.
\end{itemize}
\end{prop}

\begin{remark}
The slope of the boundary component of the annulus on $\bdry N(L)$ is taken with respect to the standard meridian and the longitude induced from the once-punctured torus fiber.  
\end{remark}

\begin{proof}
The annuli for $n(1)=1$ and $n(1)=0$ are shown in 
Figures~\ref{fig:longitude_annulus} and~\ref{fig:meridian_annulus} where the twisted saddles $C_{{\mathit a}, +1}$ and $C_{{\mathit a}, 0}$ respectively are stacked upon their predecessors.  Their slopes on $\bdry N(L)$ can be seen to be $-1/1$ and $1/0$ respectively.  Similarly the annuli for general $n(1)$ and their slopes of $-1/n(1)$ on $\bdry N(L)$ are obtained.

\begin{figure}
\centering
\input{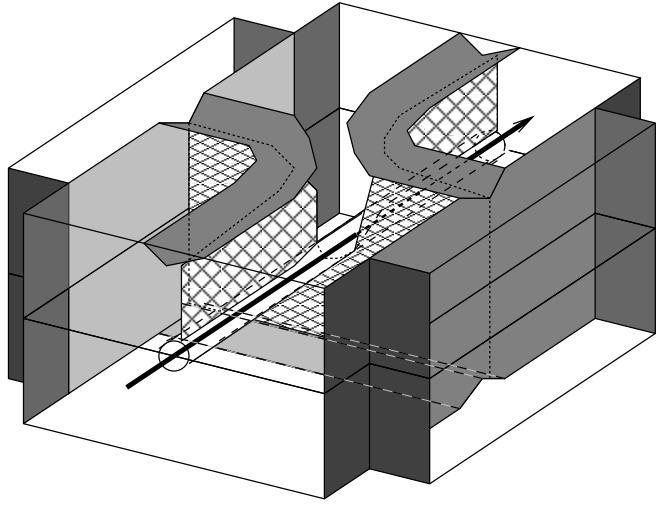}
\caption{Longitudinal annulus in $C(0; n(k), \dots, n(2), +1)$.}
\label{fig:longitude_annulus}
\end{figure}

\begin{figure}
\centering
\input{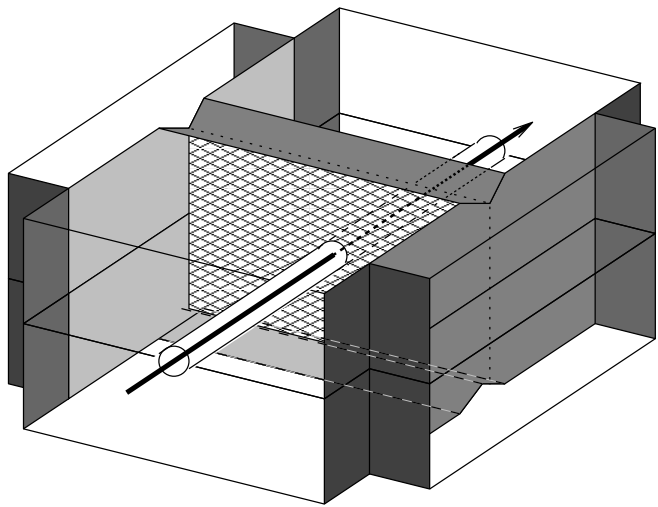}
\caption{Meridional annulus in $C(0; n(k), \dots, n(2), 0)$.}
\label{fig:meridian_annulus}
\end{figure}

By \cite[Remark 2.5.2]{cjr:isioptb}, $\widehat{C}(0; n(k), \dots, n(2), n(1))$ is incompressible in $\widehat{M}$ if $|n(1)| \geq 2$.
That $\widehat{C}(0; n(k), \dots, n(2), \pm 1)$ bounds a handlebody in $\widehat{M}$ whereas $\widehat{C}(0; n(k), \dots, n(2), \pm 1)$ does not is, however, harder to see.  

{\em Case 1.}
If the incompressible surface $\widehat{C}(0; n(k), \dots, n(2), n(1))$ in $\widehat{M}-N(L)$ were to bound a handlebody in $\widehat{M}$, then $L$ must be contained in the handlebody.  Therefore we will focus on the component containing $L$.

Consider the essential surface $C = C(0; n(k), \dots, n(2), \pm 1)$ with meridional boundary in $M-N(L)$ and the corresponding surface $\widehat{C} = \widehat{C}(0; n(k), \dots, n(2), \pm 1)$ in $\widehat{M}-N(L)$ that is disjoint from $L$ but intersects $K$ four times.  Filling the boundary component of $M-N(L)$ corresponding to $L$ trivially (or ``forgetting'' $L$) makes $C$ boundary compressible in $M \subset \widehat{M}$.  After the first two boundary compressions indicated in Figure~\ref{fig:boundary_compressions}, the boundary of the resulting surface has two components that bound disks on $\bdry M$.  Capping off these two boundary components with disks yields a closed surface disjoint from $K$, see Figure~\ref{fig:cappedends}.  Indeed, this surface is isotopic in $\widehat{M}$ to $\widehat{C}$.   The two compressing disks guide isotopies of $\widehat{C}$ through $K$ (and $L$).
\begin{figure}
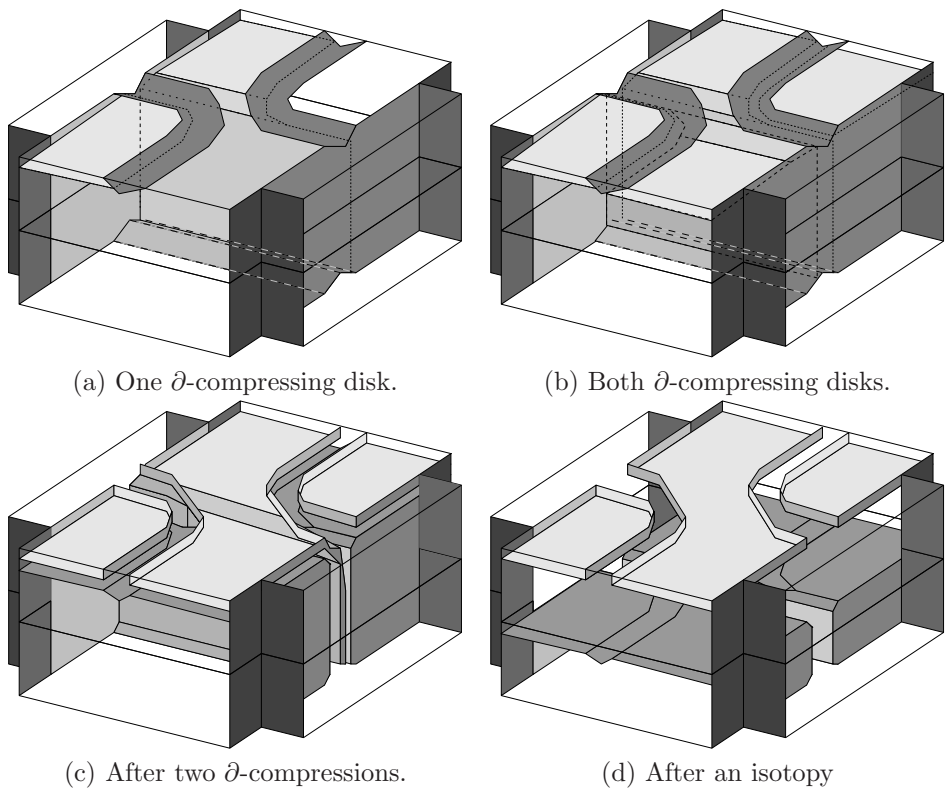

\centering
$\begin{array}{cc}
\input{Figures/oneboundarycompression.pstex_t} & \input{Figures/twoboundarycompressions.pstex_t}\\
\mbox{(a) One $\bdry$-compressing disk.} & \mbox{(b) Both $\bdry$-compressing disks.}\\
\input{Figures/aftercompression.pstex_t} & \input{Figures/afterisotopy.pstex_t}\\
\mbox{(c) After two $\bdry$-compressions.} & \mbox{(d) After an isotopy}
\end{array}$
\caption{Boundary compressions of $C(0; n(k), \dots, n(2), \pm 1)$.}
\label{fig:boundary_compressions}
\end{figure}
\begin{figure}
\centering
\input{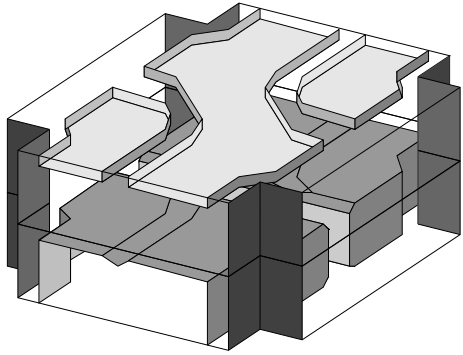}
\caption{Adding sides to the first and last block.}
\label{fig:cappedends}
\end{figure}
Viewed in the blocks of 
\[M = T \times I/_{\alpha^{\pm 1}} T \times I/_{\beta^{n(2)}} \dots T \times I/_{\beta^{n(k)}},\]
this effectively adds to each of the middle $k-2$ twisted saddles two vertical disks each parallel into $\bdry T \times I$; one with vertical edges $\bdry{\mathit a}_+ \times I$ and the other with vertical edges $\bdry{\mathit a}_- \times I$.  See Figures~\ref{fig:side_saddles} (a) and (b).
\begin{figure}
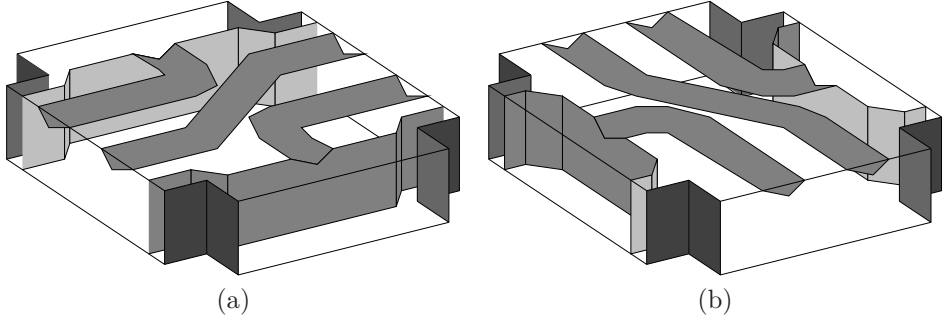

\centering
$\begin{array}{cc}
\input{Figures/side_saddle_a.pstex_t} &
\input{Figures/side_saddle_b.pstex_t} \\
\mbox{(a)} &
\mbox{(b)} 
\end{array}$
\caption{Adding sides to twisted saddles.}
\label{fig:side_saddles}
\end{figure}
After this isotopy of the surface through $K$, we can see it bounding a handlebody as it sits in $M$.  Now that the twisted saddles have sides, there is a compressing disk on the fiber joining the top of even blocks to the bottom of odd blocks, i.e.\ the fiber $T \times \{1\}$ of the $i$th block for $2 \leq i \leq k-2$ even.  This compressing disk is visible in $T\times\{1\}$ of Figure~\ref{fig:side_saddles} (a) and in $T \times \{0\}$ of Figure~\ref{fig:side_saddles} (b).

Compressing the surface along these $\frac{1}{2}k-1$ disks separates the surface into $\frac{1}{2}k$ tori.  There is one torus in each pair of an odd numbered block with its successive even numbered block.  These tori bound solid tori.  Undoing the $\frac{1}{2}k-1$ compressions by attaching $1$--handles to these solid tori and reversing the isotopies, we conclude that $\widehat{C}=\widehat{C}(0;n(k), \dots, n(2), \pm 1)$ bounds a handlebody.

{\em Case 2.}
Consider the surface $\widehat{C} = \widehat{C}(0;n(k), \dots, n(2), 0) \subset \widehat{M}$ and the corresponding surface $C=C(0;n(k), \dots, n(2), 0) \subset M$ in the complement of $L$.  Let $H$ be the $\SL_2(\Z)$ monodromy of the punctured torus bundle $M$.  Thus $H$ is conjugate to $B^{n(k)} A^{n(k-1)} \dots A^{n(3)} B^{n(2)} A^0$.  Therefore
\begin{align*}
[H] &= [ B^{n(k)} A^{n(k-1)} \dots A^{n(3)} B^{n(2)} A^0 ] \\
    &= [ A^{n(k-1)} \dots A^{n(3)} B^{n(2)} A^0 B^{n(k)} ] \\
    &= [ A^{n(k-1)} \dots A^{n(3)} B^{n(2) + n(k)} ] \\
    &= [ P^{-1} A^{n(k-1)} \dots A^{n(3)} B^{n(2) + n(k)} P ] \\
    &= [ P^{-1} P B^{n(k-1)} \dots B^{n(3)} A^{n(2) + n(k)} ]\\
    &= [ B^{n(k-1)} \dots B^{n(3)} A^{n(2) + n(k)} ].
\end{align*}
Since $|n(i)| \geq 2$ for $i = 3, \dots, k-1$ and the exponent sum is no different, by Propositions~\ref{prop:2.5.1} and~\ref{prop:meridian_conditions} this corresponds to another essential surface $C(0;n(k-1), \dots, n(3), n(2)+n(k))$ with meridional boundary in $M$ in the complement of another level knot $L'$.  This surface $C(0;n(k-1), \dots, n(3), n(2)+n(k))$ is obtained by compressing along the ``obvious'' compressing disk $D$ for $C(0;n(k), \dots, n(2), 0)$ forgetting $L$.  Note that $L'$ must be contained in the component of $M - N(C(0;n(k-1), \dots, n(3), n(2)+n(k)))$ that does not come from the component of $M - N(C(0;n(k), \dots, n(2), 0))$ that contained $L$.  

By Theorem~\ref{thm:2.5.2}, $\widehat{C}(0;n(k-1), \dots, n(3), n(2)+n(k))$ is incompressible in $\widehat{M}-N(L)$.
Also note that by Lemma~\ref{lem:n1leq1}, $n(2)+n(k) = \pm 1$ or $0$ if $\widehat{C}(0;n(k-1), \dots, n(3), n(2)+n(k))$ is compressible in $\widehat{M}$.

Assume $\widehat{C}(0;n(k), \dots, n(2), 0)$ bounds a handlebody $V$ to the side containing $L$ in $\widehat{M}$.  The other side must be incompressible.  Note that $\widehat{C}(0;n(k), \dots, n(2), 0) \not \cong S^2$ since $D$ is a properly embedded separating disk that $L$ intersects once.  The compression along $D$ yielding $\widehat{C}(0;n(k-1), \dots, n(3), n(2)+n(k))$ compresses $V$ along a nonseparating curve to yield another handlebody $V' \subset V$ which does not contain the new level knot $L'$.  This contradicts that $\widehat{C}(0;n(k-1), \dots, n(3), n(2)+n(k))$ is essential in $\widehat{M}-N(L')$ unless $\widehat{C}(0;n(k-1), \dots, n(3), n(2)+n(k)) \cong S^2$.  

But if this resulting surface is $S^2$, then $V' \cong B^3$ and $V$ is a solid torus.  Thus $\widehat{C}(0;n(k), \dots, n(2), 0)$ is parallel to $\bdry N(L)$ and is not essential.
\end{proof}

Indeed if $\widehat{M} \cong S^3$, as one may check, there are no special forms $B^{n'(2)}A^{n'(1)}$ or $B^{n'(4)}A^{n'(3)}B^{n'(2)}A^{n'(1)}$ where $|n'(i)| \geq 2$ for $i \geq 2$ and $|n'(1)| \leq 1$ that represent the monodromy of a genus one fibered knot.


\section{The algorithm for Berge knots}\label{sec:bergeknotalg}

In preparation for the proof of Theorem~\ref{thm:generalapplication} we streamline Algorithms~\ref{algor:main_algorithm} and Algorithm~\ref{algor:framing_algorithm} for the Berge knots that lie on the fiber of a trefoil or the figure eight knot.

Let $K$ be the left handed trefoil, right handed trefoil, or figure eight knot.  Then $S^3-N(K) \cong T \times I /_\eta$ where $\eta$ is $\alpha^{-1} \circ \beta^{-1}$, $\beta \circ \alpha$, or $\alpha \circ \beta^{-1}$ respectively.
Let $H \in \SL_2(\Z)$ be the matrix associated to $\eta$ so that $H$ is $A^{-1}B^{-1}$, $BA$, or $AB^{-1}$ respectively. The corresponding standard meridional framing is $\eta$.
 
\begin{algor}\label{algor:streamlined_alg}
Let $K$ and $\eta$ be as above.  Let $L \in \K_\eta$ be given by the slope $x/y$ with, say, $y \geq 0$.  The following steps give a procedure to list all closed connected essential surfaces in $S^3-N(L)$. 

\begin{enumerate}
\item[{\em Step 1.}]  Choose a change of basis matrix $W \in \SL_2(\Z)$ so that $W \begin{pmatrix} 1 \\ 0 \end{pmatrix} = \begin{pmatrix} x \\ y \end{pmatrix}$.

\item[{\em Step 2.}]  Let $X(N) = (W A^N)^{-1} H (W A^N)$.

\item[{\em Step 3.}]  List all $N$ such that $X(N) = \begin{pmatrix} p(N) & r(N) \\ q(N) & s(N) \end{pmatrix}$ where $|p(N)| < |q(N)|$.

\item[{\em Step 4.}]  For each $N$ listed in Step 3,
\begin{itemize} 
\item obtain the MCFEs $[a_1, a_2, \dots, a_{k-1}]$ for $p(N)/q(N)$ such that $k$ is even, and
\item find $a_k$ such that $[a_1, a_2, \dots, a_{k-1}, a_k]$ is a continued fraction expansion for $r(N)/s(N)$.
\end{itemize}

\item[{\em Step 5.}]  For each $N$ of Step 3 and for each minimal continued fraction expansion $[a_1, a_2, \dots, a_{k-1}, a_k]$ obtained in Step 4, list those such that
\[ \sum a_i = E(H). \]

\item[{\em Step 6.}]  For each continued fraction $[a_1, a_2, \dots, a_{k-1}, a_k]$ listed in Step 5 there is a closed essential surface of type $\widehat{C}(0; a_1, a_2, \dots, a_{k-1}, a_k)$ in $S^3-N(L)$.
\end{enumerate}
\end{algor}

\begin{proof}
Assume $\widehat{S}$ is a closed essential surface in $S^3_L$.  Isotop $\widehat{S}$ to intersect $K$ minimally.  Let $S = \widehat{S} - N(K)$.  By Lemma~\ref{lem:ess_with_bdry} $S$ is essential in $S^3_{K \cup L}$.

Since $S$ has meridional boundary on $\bdry N(K)$, $\bdry S$ is transverse to the induced fibration on $\bdry N(K)$.  From this and that $\eta$ (for each of the three choices) does not fix the isotopy class of any essential curve on $T$, Theorem~\ref{thm:surface_type} implies $S$ is isotopic to a surface of type $C(J; n(k), \dots, n(1))$.  Since $S$ is essential, Proposition~\ref{prop:2.5.1} implies that $|n(i)| \geq 2$ for $i = 2, \dots, k$.  Since the components of $\bdry S$ are meridional curves, $J = 0$ and $\sum_{i=1}^{k} n(i) = E(H)$ by Proposition~\ref{prop:meridian_conditions}.

Therefore, to every closed essential surface in $S^3_L$ there is a surface of type $C(0; n(k), \dots, n(1))$ in $S^3_{K \cup L}$ such that $\sum_{i=1}^{k} n(i) = E(H)$.  Algorithm~\ref{algor:main_algorithm} lists all surfaces of type $C(J; n(k), \dots, n(1))$ that are essential in $S^3_{K \cup L}$.  The meridional conditions $J = 0$ and $\sum_{i=1}^{k} n(i) = E(H)$ are then simple to check.


We now explain the shortcuts beginning at Step 4 where the algorithm first diverges from Algorithm~\ref{algor:main_algorithm}.

Since $S$ must have the type of an essential twisted surface with meridional boundary components, $J=0$.  Hence in Step 4 we take the continued fraction expansions $[a_1, a_2, \dots, a_{k-1}]$ where $k$ is even.  

At the end of Step 4, we may conclude that $X(N) = P^J B^{a_1} \dots B^{a_{k-1}} A^{a_k}$ where $J = 0$ or $J = 2$.  If $J = 2$, this special form has standard framing $\xi = \phi^2 \beta^{a_1} \dots \beta^{a_{k-1}} \alpha^{a_k}$ and since $X(N)$ is conjugate to $H$ which has framing $\eta$, Lemma~\ref{lem:framing} states that
\[ (E_{\phi}(\xi \eta^{-1}),E_{\psi}(\xi \eta^{-1})) = (4j, 0) \mod (2,3).\]
But since Step 5 requires $\sum a_i = E(H)$,
 \[E_\phi(\eta) = E_\phi(\xi)-2 \mbox{ and } E_\psi(\eta) =E _\psi(\xi).\]
Thus
 \[(E_{\phi}(\xi \eta^{-1}),E_{\psi}(\xi \eta^{-1})) = (-2, 0) \mod (2,3),\]
and we get a contradiction.  Therefore $J=0$.

As in Step 5 of Algorithm~\ref{algor:main_algorithm}, having a special form $B^{a_1} \dots B^{a_{k-1}} A^{a_k}$ for $X(N)$ implies the existence of an essential surface of type $C(0;a_1, \dots, a_k)$.  The conditions of the current Step 5 imply that the boundary of this surface is meridional.  

As stated in Step 6, the surface then caps off to a closed surface $\widehat{C}(0;a_1, \dots, a_k)$ under the meridional filling of the boundary component of $S^3_{K \cup L}$ corresponding to $K$.  By Theorem~\ref{thm:2.5.2} this surface is essential.
\end{proof}

\subsection{Examples}

Let $K$ be the left handed trefoil.  Let $\eta = \alpha^{-1} \circ \beta^{-1}$ and $H=A^{-1}B^{-1}$.  
Here we give some examples of the application of Algorithm~\ref{algor:streamlined_alg} to knots in $\K_\eta = \K_{\alpha^{-1} \circ \beta^{-1}}$, the knots on the fiber of the left handed trefoil.

\begin{example}\label{exa:genus2sfces}
Let $L \in \K_\eta$ be a knot with slope $\frac{1+2z}{4z}$ such that $z$ is an integer $\geq 3$.   We apply Algorithm~\ref{algor:streamlined_alg}.
\begin{enumerate}
\item[{\em Step 1.}] Choose $W = \begin{pmatrix} 1+2z & -z \\ 4z & 1-2z \end{pmatrix}$.

\item[{\em Step 2.}] $X(N) = (W A^N)^{-1} H (W A^N)$.
  {\footnotesize   \[X(N)=
  \begin{pmatrix} 3z-6z^2 + (1 + 12z^2)N & 1-3z+3z^2 -(1-6z+12z^2)N +(1 + 12z^2)N^2 \\
  -(1 + 12z^2) & 1-3z+6z^2 - (1 + 12z^2)N \end{pmatrix} \]
}
\item[{\em Step 3.}] $q(N) = -(1+12z^2)$ and $p(N) = 3z-6z^2 + (1 + 12z^2)N$.  Therefore $|p(N)| < |q(N)|$ if $N = 0$ or $N = -1$.

\item[{\em Step 4.}]  Recall that we are assuming $z\geq 3$. 
\begin{enumerate}
  \item[($N=0$.)] $p(0)/q(0) = \frac{-3z+6z^2}{1 + 12z^2}$ has SCFE $[2,1-z,2,-1,z-1,-2]$.  We have the following MCFEs and determine their corresponding $a_k$:  
\begin{itemize}
   \item $[2,-z,-2,-2,z,2]$, $a_k=+1$
   \item $[2,-z,-2,-3,\underbrace{-2, \dots,-2}_{z-2},-3]$, $a_k=0$
   \item $[2,1-z,3,z,-2]$, $a_k=0$
   \item $[2,1-z,3,z+1,2]$, $a_k=+1$
   \item $[3,\underbrace{2, \dots, 2}_{z-2},3,\underbrace{-2, \dots,-2}_{z-1},-3]$, $a_k=0$
   \item $[3,\underbrace{2, \dots, 2}_{z-2},4,z,-2]$, $a_k=0$
   \item $[3,\underbrace{2, \dots, 2}_{z-2},4,z+1,2]$, $a_k=+1$
\end{itemize}

 \item[($N=-1$.)]$p(-1)/q(-1) = \frac{-1-3z-6z^2}{1 + 12z^2}$ has SCFE $[-1, 1,1-z,2,-1,z-1,-2]$.  We have the following MCFEs and determine their corresponding $a_k$: 
\begin{itemize}
   \item $[-2,-z-1,-2,-3,\underbrace{-2, \dots,-2}_{z-2},-3]$, $a_k=-1$
   \item $[-2,-z-1,-2,-2,z,2]$, $a_k=0$
   \item $[-2,-z-1,-2,-2,z-1,-2]$, $a_k=-1$
   \item $[-2,-z,2,\underbrace{-2, \dots,-2}_{z-1},-3]$, $a_k=-1$
   \item $[-2,-z,3,z+1,2]$, $a_k=0$
   \item $[-2,-z,3,z,-2]$, $a_k=-1$
\end{itemize}
Note that some of these MCFEs have $k$ odd.  
 
\end{enumerate}

\item[{\em Step 5.}] Since $E(H)=E(A^{-1}B^{-1})=-2$, we obtain just one list $\{a_1, \dots, a_{k-1},a_k\}$ such that $k$ is even and $\sum_{i=1}^k a_i = -2$:
\[ \{-2,-z,3,z,-2,-1\} \]

\item[{\em Step 6.}] The knot $L$ given by $(1+2z, 4z)$ for $z\geq 3$ contains the closed essential surface $\widehat{C}(0;-2,-z,3,z,-2,-1)$ in its complement. Furthermore, this is the only one.

\end{enumerate}

\end{example}

\begin{remark}
One can actually show that a knot in $\K_\eta$ with slope $\frac{1+2z}{4z}$ for $|z| \geq 2$ has the surface $\widehat{C}(0;-2,-z,3,z,-2,-1)$ as the only closed essential surface in its complement, but the restriction $z \geq 3$ eliminates special cases that would arise otherwise.  Moreover, one can show that if a knot in $\K_\eta$ has a closed essential genus $2$ surface in its complement, then it is one of these knots.  A similar result is also true for the knots on the fiber of the figure eight knot.
\end{remark}

\begin{example}\label{exa:highergenus}

By doing the reverse of the compression in Case 2, Proposition~\ref{prop:differences_in_surfaces} we can construct many more knots with close essential surfaces in their complements.  We use the above example as a model for the general construction.   In Example~\ref{exa:genus2sfces} we obtained the closed essential surface $\widehat{C}(0;-2,-z,3,z,-2,-1)$ in the complement of the knot $\frac{1+2z}{4z} \in \K_\eta$.  Thus
{\small \begin{align*} 
(A^{-N} W^{-1}) A^{-1}B^{-1} (W A^N) &= B^{-2} A^{-z} B^{3} A^{z} B^{-2} A^{-1} \\
                          &= B^{-2} A^{-z} B^{3} A^{z} B^{-2} A^{z'-1} B^0 A^{-z'}\\
A^{-z'}(A^{-N} W^{-1}) A^{-1}B^{-1} (W A^N)A^{z'} &= A^{-z'} B^{-2} A^{-z} B^{3} A^z B^{-2} A^{z'-1} B^0 \\
P^{-1}A^{-z'}(A^{-N} W^{-1}) A^{-1}B^{-1} (W A^N)A^{z'}P &= P^{-1}A^{-z'} B^{-2} A^{-z} B^{3} A^z B^{-2} A^{z'-1} B^0 P \\
(P^{-1}A^{-z'}A^{-N} W^{-1}) A^{-1}B^{-1} (W A^N A^{z'}P) &= B^{-z'} A^{-2} B^{-z} A^{3} B^z A^{-2} B^{z'-1} A^0. 
\end{align*} }
Since $W = \begin{pmatrix} 1+2z & -z \\ 4z & 1-2z \end{pmatrix}$ and $N = 0$, 
\[W A^N A^{z'}P = \begin{pmatrix} -z-z'-2z z' & -1-2z \\
                                 1-2z-4z z' & -4z \end{pmatrix}.\]
Hence if $|z'| \geq 2$ and $|z'-1| \geq 2$ (and $|z| \geq 2$), the knot $\frac{z+z'+2z z'}{-1+2z+4z z'} \in \K_\eta$ has in its complement the closed essential surface $\widehat{C}(0;-z', -2, -z, 3, z, -2, z'-1, 0)$ of genus $3$. 

This process may be repeated to construct knots in $\K_\eta$ whose complements contain a closed essential surface of genus $g$ for any $g \geq 2$.  Where we had replaced $A^{-1}$ with $A^{z'-1} B^{0} A^{-z'}$, one would now replace $A^0$ with $A^{z''} B^0 A^{-z''}$.
\end{example}

\begin{remark}
Every surface constructed in this manner has $n(1)= 0$.  There are closed essential surfaces of every genus greater than $2$ in the complements of knots in $\K_\eta$, however, that do not arise this way.  
\end{remark}

\begin{example}\label{exa:smallknot}

As we shall see in Lemma~\ref{lem:smallknots}, there are many small knots in $\K_\eta$.  Here we give an example of the algorithm applied to a simple knot (which is not a torus knot) to show it is small.

Let $L \in \K_\eta$ be the knot $-\frac{3}{2} \in \K_\eta$.  We apply Algorithm~\ref{algor:streamlined_alg}.
\begin{enumerate}
\item[{\em Step 1.}] Choose $W = \begin{pmatrix} -3 & 1 \\ 2 & -1 \end{pmatrix}$.

\item[{\em Step 2.}] $X(N) = (W A^N)^{-1} H (W A^N)$.
     \[ X(N) =
 \begin{pmatrix} -7 - 19N & 3+15N+19N^2 \\
                 -19      & 8+19N        \end{pmatrix} \]

\item[{\em Step 3.}] $q(N) = -19$ and $p(N) = -7 - 19N$.  Therefore $|p(N)| < |q(N)|$ if $N = 0$ or $N = -1$.

\item[{\em Step 4.}] $\mbox{}$
\begin{enumerate}
  \item[($N=0$.)] $p(0)/q(0) = \frac{7}{19}$ has SCFE $[2,-1,2,-2]$ and two MCFEs of odd length:  $[3,4,2]$ and $[3,3,-2]$.  

For $[3,4,2]$:
\begin{align*}
 [3,4]      &= (r(0)+ p(0) a_4)/(s(0)+ q(0) a_4) \\
\Rightarrow \frac{4}{13} &= \frac{3-7a_4}{8-19a_4}\\
 \Rightarrow a_4 &= 1
\end{align*}
For $[3,3,-2]$
\begin{align*}
 [3,3]      &= (r(0)+ p(0) a_4)/(s(0)+ q(0) a_4) \\
\Rightarrow \frac{3}{8} &= \frac{3-7a_4}{8-19a_4}\\
 \Rightarrow a_4 &= 0
\end{align*}

  \item[($N=-1.$)]$p(-1)/q(-1) = -\frac{12}{19}$ has SCFE $[-1,1,-1,2,-2]$ and no MCFE of odd length.
\end{enumerate}

\item[{\em Step 5.}] We have two lists: 
\[ \{3,4,2,1\} \mbox{ and } \{3,3,-2,0\}. \] 
Their sums are $10$ and $4$ respectively.  Neither equals $E(H)=-2$.  

\item[{\em Step 6.}] Since we have no lists satisfying all the criteria, there are no closed essential surfaces in the complement of $L$.
\end{enumerate}
\end{example}



\section{The general equation} \label{sec:generalities}
We apply Algorithm~\ref{algor:streamlined_alg} to the knots $L \in \K_\eta$ in general.  Throughout this section when $N = \pm1$, the symbol $\pm$ agrees with the sign of $N$. 

Let $[\bit{b}] = [b_1, b_2, \dots, b_k]$ be the SCFE for the slope $\frac{x}{y}$ of the knot $L$.  Recall that only $b_1$ may be $0$ and consecutive coefficients do not have the same sign.  We also assume $\frac{x}{y} \neq 0$, $\pm 1$, or $\frac{1}{0}$ since these slopes all correspond to the unknot.

\begin{thm}\label{thm:generalapplication}
For $L$ as above, every closed essential surface in its complement corresponds to a solution of one of the following equations:
\begin{enumerate}
\item If $b_1 \neq 0$ or $1$,
\[ 0 =
   \sum_{i \in I} -b_i + 
   \sum_{j \in J} b_j +
   \begin{cases}
      0 & \mbox{ if } 1 \in I \\
      -1 & \mbox{ otherwise}
   \end{cases}
\]
where $I$ and $J$ are subsets of $\{1, \dots, k\}$ each not containing consecutive integers and $1 \not \in I \cap J$. 

Furthermore: 
\begin{itemize}
\item If $b_l = \pm 1$ for $l \geq 2$, then $\{l-1, l, l+1\} \cap I \neq \emptyset$ and  $\{l-1, l, l+1\} \cap J \neq \emptyset$.  
\item If $b_1 = -1$, then either $1 \in I$ or $\{1,2\} \cap J \neq \emptyset$.  
\item If $b_1 = 2$, then either $1 \in J$ or $\{1,2\} \cap I \neq \emptyset$.
\end{itemize}
\item If $b_1 = 0$ and $b_2 \neq -1$,
\[ 0 =
   \sum_{i \in I} -b_i + 
   \sum_{j \in J} b_j +
   \begin{cases}
      0 & \mbox{ if } 2 \in J \\
      -1 & \mbox{ otherwise}
   \end{cases}
\]
where $I$ and $J$ are subsets of $\{2, \dots, k\}$ each not containing consecutive integers and $2 \not \in I \cap J$.

Furthermore:
\begin{itemize}
\item If $b_l = \pm 1$ for $l \geq 3$, then $\{l-1, l, l+1\} \cap I \neq \emptyset$ and  $\{l-1, l, l+1\} \cap J \neq \emptyset$.
\item If $b_2 = -2$ then either $2 \in I$ or $\{2,3\} \cap J \neq \emptyset$.  
\item If $b_2 = 1$ then either $2 \in J$ or $\{2,3\} \cap I \neq \emptyset$.
\end{itemize}
\item If $b_1 = 0$ and $b_2 = -1$,
\[ 0 =
   \sum_{i \in I} -b_i + 
   \sum_{j \in J} b_j +
   \begin{cases}
      0 & \mbox{ if } 3 \in J \\
      -1 & \mbox{ otherwise}
   \end{cases}
\]
where $I$ and $J$ are subsets of $\{3, \dots, k\}$ each not containing consecutive integers and $3 \not \in I \cap J$. 

Furthermore:
\begin{itemize}
\item If $b_l = \pm 1$ for $l \geq 4$, then $\{l-1, l, l+1\} \cap I \neq \emptyset$ and  $\{l-1, l, l+1\} \cap J \neq \emptyset$. 
\item If $b_3 = 1$ then either $3 \in J$ or $\{3,4\} \cap I \neq \emptyset$. 
\end{itemize}
\item If $b_1 = 1$,
\[ 0 =
   \sum_{i \in I} -b_i + 
   \sum_{j \in J} b_j +
   \begin{cases}
      0 & \mbox{ if } 2 \in I \\
      -1 & \mbox{ otherwise}
   \end{cases}
\]
where $I$ and $J$ are subsets of $\{2, \dots, k\}$ each not containing consecutive integers and $2 \not \in I \cap J$.

Furthermore:
\begin{itemize}
\item If $b_l = \pm 1$ for $l \geq 3$, then $\{l-1, l, l+1\} \cap I \neq \emptyset$ and  $\{l-1, l, l+1\} \cap J \neq \emptyset$.
\item If $b_2 = -1$ then either $2 \in I$ or $\{2,3\} \cap J \neq \emptyset$.
\end{itemize}
\end{enumerate}
\end{thm}

\begin{remark}
If $[\bit{b}]$ is a MCFE as well as a SCFE with $b_1 \neq 2$ then Cases (2), (3), and (4) and the extra conditions on the sets $I$ and $J$ in Case (1) are not needed.
\end{remark}

\begin{proof}
Set $W = \mat{x}{y}{t}{u} \in \SL_2(\Z)$ where
\begin{align*}
\frac{x}{y} &= [b_1, b_2, \dots, b_k] \mbox{ and} \\
\frac{t}{u} &= [b_1, b_2, \dots, b_{k-1}].
\end{align*}
By Lemma~\ref{lem:matrix_cfe1} $W$ may then be written as 
\[ \pm B^{b_1} A^{b_2} \dots B^{b_k} \mbox{ or } \pm B A B A^{b_1} B^{b_2} \dots B^{b_k} \]
depending on the parity of $k$. 

Let
\[X(N) = (W A^N)^{-1} A^{-1}B^{-1} (W A^N) = \mat{p(N)}{r(N)}{q(N)}{s(N)} A^N\] 
for $N \in \Z$.  Then $X(N)$ may be written as either
\[
\begin{split}
(\pm A^{-N} B^{-b_k} \dots A^{-b_2} B^{-b_1}) A^{-1} B^{-1} (\pm B^{b_1} A^{b_2} \dots B^{b_k} A^N) \\
= A^{-N} B^{-b_k} \dots A^{-b_2} B^{-b_1} A^{-1} B^{b_1 - 1} A^{b_2} \dots B^{b_k} A^N
\end{split}
\]
or 
\[
\begin{split}
(\pm A^{-N} B^{-b_k} \dots B^{-b_2} A^{-b_1} B^{-1} A^{-1} B^{-1}) A^{-1} B^{-1} (\pm B^1 A^1 B^1 A^{b_1} B^{b_2} \dots B^{b_k} A^N) \\
= A^{-N} B^{-b_k} \dots B^{-b_2} A^{-b_1} B^{-1} A^{b_1 - 1} B^{b_2} \dots B^{b_k} A^N
\end{split}
\]
depending on the parity of $k$.  In either case, note that 
\[ p(N)/q(N) = [0, -N, -b_k, \dots, -b_2, -b_1, -1, b_1 -1 , b_2, \dots, b_k]. \]
Define $[\bit{x}(N)] = [0, -N, -b_k, \dots, -b_2, -b_1, -1, b_1 -1 , b_2, \dots, b_k]$.

According to Algorithm~\ref{algor:streamlined_alg} for every MCFE $[\bit{x}_m] = [a_1, a_2, \dots, a_l]$ of $p(N)/q(N)$ of odd length such that $N' = -2 - \sigma(\bit{x}_m)$ and $[a_1, a_2, \dots, a_l, N']$ is a continued fraction expansion for $r(N)/s(N)$ for some $N$ there is a closed essential surface in the complement of $L$.  Since every MCFE of a rational number is derived from its SCFE, we need to obtain the SCFE $[\bit{x}_s(N)]$ for $p(N)/q(N)$.  We determine the SCFE in \S\ref{subsec:SCFEforp/q} and the MCFEs corresponding to essential surfaces in \S\ref{subsec:MCFEsforp/q}.

\subsection{The SCFE for $p(N)/q(N)$.} \label{subsec:SCFEforp/q}

Recall that if $p(N)/q(N)$ is to have a MCFE then its absolute value must be less than $1$.  Since $p(N)/q(N) = [\bit{x}(N)]$,
\[
 p(N)/q(N) = \cfrac{1}{0 - \cfrac{1}{-N - [\bit{b'}]}} 
           = N + [\bit{b'}]
\]
where $[\bit{b'}] = [-b_k, \dots, -b_2, -b_1, -1, b_1 -1 , b_2, \dots, b_k]$.

\begin{lemma}\label{lem:sizeofb'}
$|[\bit{b'}]| < 1$ and $\sgn([\bit{b'}]) = \sgn(-b_k)$.
\end{lemma}
\begin{proof}
We exhibit the SCFE for $[\bit{b'}]$.  The first coefficient will be nonzero and hence imply the conclusion of the lemma.

By move (CF1),
\[ [\bit{b'}] \mapsto [\bit{b}_1'] \underset{\mathrm{def}}{=} [-b_k, \dots, -b_2, -b_1 +1, b_1, b_2, \dots, b_k]. \]
This is a SCFE unless $b_1 = 0$ or $1$. 

\noindent{\em Case $b_1 = 0$.}  
By move (CF1)
\begin{align*}
[\bit{b}_1'] &= [-b_k, \dots, -b_2, 1, 0, b_2, \dots, b_k] \\
                  &\mapsto [-b_k, \dots, -b_3, -b_2, b_2 +1, b_3, \dots, b_k] \underset{\mathrm{def}}{=} [\bit{b}_2']
\end{align*}
which is a SCFE unless $b_2 = -1$ or $0$.  Since $[\bit{b}]$ is a SCFE, $b_2$ cannot be $0$.
Furthermore $k > 1$ since $[\bit{b}] \neq 0$.

If $b_2 = -1$ then by move (CF1)
\begin{align*}
[\bit{b}_2'] &= [-b_k, \dots, -b_3, 1, 0, b_3, \dots, b_k] \\
                  &\mapsto [-b_k, \dots, -b_4, -b_3, b_3 +1, b_4, \dots, b_k] \underset{\mathrm{def}}{=} [\bit{b}_3']
\end{align*}
which is a SCFE unless $b_3 = -1$ or $0$.  Since $[\bit{b}]$ is a SCFE and $b_2 = -1$ neither of these two situations can occur.  Furthermore if $k = 2$ then $[\bit{b}] = [0, -1] = 1$ contrary to our assumption on $\frac{x}{y}$.

\noindent{\em Case $b_1 = 1$.}
By move (CF1)
\begin{align*}
[\bit{b}_1'] &= [-b_k, \dots, -b_2, 0, 1, b_2, \dots, b_k] \\
                  &\mapsto [-b_k, \dots, -b_3, -b_2 +1, b_2, b_3, \dots, b_k] \underset{\mathrm{def}}{=} [\bit{b}_4']
\end{align*}
which is a SCFE unless $b_2 = 0$ or $1$.  Since $[\bit{b}]$ is a SCFE and $b_1 = 1$ neither of these two situations can occur.  Furthermore if $k = 1$ then $[\bit{b}] = 1$ contrary to our assumption on $\frac{x}{y}$.

Each $[\bit{b}_i']$, for $i = 1, 2, 3, 4$, has nonzero leading coefficient.  One of these is the SCFE for $[\bit{b'}]$.  Therefore by Lemma~\ref{lem:scfesize}, $|[\bit{b'}]|<1$ and $\sgn([\bit{b'}]) = \sgn(-b_k)$.
\end{proof}

\begin{cor}\label{cor:whichN}
If $p(N)/q(N)$ has a MCFE then either $N = 0$ or $N = \sgn(b_k)$.
\end{cor}
\begin{proof}
A rational number has a MCFE only if its absolute value is less than $1$.
Since $p(N)/q(N) = [0,N,\bit{b'}] = N + [\bit{b'}]$ and by Lemma~\ref{lem:sizeofb'} $|[\bit{b'}]| < 1$, $|p(N)/q(N)| < 1$ if and only if either $N = 0$ or $N = \pm 1$ such that $\sgn(N) \neq \sgn([\bit{b'}])$.   However, $\sgn([\bit{b'}]) = \sgn(-b_k) = -\sgn(b_k)$. 
\end{proof}

The above corollary shows that we really only need to consider two values of $N$ in Algorithm~\ref{algor:streamlined_alg}.

We now find the SCFE $[\bit{x}_s(N)]$ for $p(N)/q(N)$.  

If $N=0$, then for one of $i = 1, 2, 3,$ or $4$.
\[p(0)/q(0) = [\bit{x}(0)]= [0, 0, \bit{b}_i'] = [\bit{b}_i']\]
is the SCFE.

If $N=\pm 1$, then for $i = 1, 2, 3,$ or $4$ such that $[\bit{b}_i']$ is a SCFE
\begin{align*}
 p(\pm1)/q(\pm1) = [\bit{x}(\pm1)]=[0, \mp 1, \bit{b}_i'] &= [0, \mp 1, -b_k, \dots, b_k] \\
                                         &\mapsto [\pm 1, -b_k \pm 1, \dots, b_k] \underset{\mathrm{def}}{=} [\bit{b}_i'']
\end{align*}
by move (CF1).  This is a SCFE unless $b_k = 0$ or $\pm 1$.  Since $[\bit{b}_i']$ is a SCFE, $|b_k| \geq 2$.  Therefore $[\bit{b}_i'']$ is the SCFE.

Notice that the move (CF1') is never used in obtaining $[\bit{x}_s(N)]$ from $[\bit{x}(N)]$.  Hence their last two partial sums are equal.  Furthermore, as a consequence of Lemma~\ref{lem:matrix_cfe2} there exists an $N_s$ such that $[\bit{x}(N),N] = [\bit{x}_s(N),N_s]$.  By Lemma~\ref{lem:threesums}, $N_s = N'$.

There are two extra ``degenerate'' cases to consider for the SCFE of $p(N)/q(N)$.  

If $n=2$, $i=4$ (so $b_1 = 1$), and $N=1$, then 
\[ [\bit{b}_4'] = [1, -b_2 +1 +1, b_2] = [1, -b_2+2, b_2] \]
is a SCFE unless $b_2 = 0, 1,$ or $2$.  These cannot occur since $b_1 =1$ and $[\bit{b}]$ is a SCFE.

If $n=1$, $i=1$, and $N=1$, then
\[ [\bit{b}_1'] = [1, -b_1 +1 +1, b_1] = [1, -b_1 +2, b_1] \]
is a SCFE unless $b_1 = 0, 1,$ or $2$.  By our assumption on $\frac{x}{y}$, $b_1 \neq 0$ or $1$.  However, $b_1$ may equal $2$ in which case $\frac{x}{y} = \frac{1}{2}$.  This slope corresponds to the unknot.

Therefore we obtain:
\begin{lemma} \label{lem:SCFEforp/q}
For $N = 0$ or $\pm 1$, 
 $p(N)/q(N)$ has SCFE $[\bit{x}_s(N)]$ where $\bit{x}_s(0) = \bit{b}_i'$ and $\bit{x}_s(\pm1) = \bit{b}_i''$  for some $i \in \{1, 2, 3, 4\}$.  Furthermore $r(N)/s(N)$ has the continued fraction expansion $[\bit{x}_s(N),N]$.
\end{lemma}

\subsection{MCFEs and surfaces} \label{subsec:MCFEsforp/q}
Given the SCFE $[\bit{x}_s(N)]$ for $p(N)/q(N)$, we determine every MCFE $[\bit{x}_m]$ of odd length and $N'$ such that both $r(N)/s(N) = [\bit{x}_m,N']$ and $\sigma(\bit{x}_m) + N' = \sigma(\bit{x}(N)) + N$.

For such an $[\bit{x}_m]$ and $N'$, $[\bit{x}_s(N),N] = r(N)/s(N)$ by Lemma~\ref{lem:SCFEforp/q}.  Then Lemma~\ref{lem:matrix_scfemcfe} implies $\sigma(\bit{x}_m) - \sigma(\bit{x}(N)) = N - N' = \begin{cases} 0 \\ \sgn(b_k) \end{cases}$.

We have four cases for the four different possibilities of the SCFE $[\bit{x}_s(N)]$ depending on $[\bit{b}]$.  Each case has the two subcases of $N = 0$ and $N = \pm1$.  See Lemma~\ref{lem:SCFEforp/q}.
\begin{enumerate}
\item If $b_1 \neq 0,1$ then $[\bit{x}_s(0)] = [\bit{b}_1']$ and $[\bit{x}_s(\pm1)] = [\bit{b}_1'']$.
\item If $b_1 = 0$ and $b_2 \neq -1$ then $[\bit{x}_s(0)] = [\bit{b}_2']$ and $[\bit{x}_s(\pm1)] = [\bit{b}_2'']$.
\item If $b_1 = 0$ and $b_2 = -1$ then $[\bit{x}_s(0)] = [\bit{b}_3']$ and $[\bit{x}_s(\pm1)] = [\bit{b}_3'']$.
\item If $b_1 = 1$ then $[\bit{x}_s(0)] = [\bit{b}_4']$ and $[\bit{x}_s(\pm1)] = [\bit{b}_4'']$.
\end{enumerate}

\noindent{\em Case (1).} $b_1 \neq 0,1$

The SCFEs for $p(N)/q(N)$ for $N = 0, \pm1$ are  
\[ [\bit{x}_s(0)] = [\bit{b}_1'] = [-b_k, \dots, -b_2, -b_1 +1, b_1, b_2, \dots, b_k] \]
and
\[ [\bit{x}_s(\pm1)] = [\bit{b}_1''] = [\pm 1, -b_k \pm 1, \dots, -b_2, -b_1 +1, b_1, b_2, \dots, b_k] .\]

\noindent{\em Subcase $N = 0$.} 

 Assume $[\bit{x}_m]$ is a MCFE for $p(0)/q(0)$.  Then it is obtained from $[\bit{x}_s(0)] = [\bit{b}_1']$ by applying moves (M) and (M') to nonadjacent coefficients of $[\bit{b}_1']$. (Cf.\ \S\ref{subsubsec:rcfe}.)  In particular, if $\pm1$ is a coefficient of $[\bit{b}_1']$, then (M) or (M') must be applied to it or an adjacent coefficient to make the all the resulting coefficients of $[\bit{x}_m]$ not $\pm1$.

The difference of coefficient sums between $[\bit{x}_m]$ and $[\bit{b}_1']$ is 

{\footnotesize
\begin{align*}
 \sigma(\bit{x}_m) - \sigma(\bit{b}_1') 
&= 
   \sum_{i \in I}-3(-b_i) + 
   \begin{cases}
   -3 &\mbox{ if } 1 \in I\\
   0 &\mbox{ otherwise}
   \end{cases}         
   +
   \sum_{j \in J} -3(b_j)                             
\shoveright{ +
   \begin{cases}
   \sgn(b_k) &\mbox{ if } k \in J\\
   0 &\mbox{ otherwise}
   \end{cases} }                                        \\
&=
   -3 \left(\sum_{i \in I} -b_i + 
   \begin{cases}
   1 &\mbox{ if } 1 \in I\\
   0 &\mbox{ otherwise}
   \end{cases}  + \sum_{j \in J} b_j 
   \right)                                          
+
   \begin{cases}
   \sgn(b_k) &\mbox{ if } k \in J\\
   0 &\mbox{ otherwise}
   \end{cases}
\end{align*} }
where  $I \subset \{1, \dots, k-1\}$ and $J \subset \{1, \dots, k\}$ each not containing consecutive integers and $1 \not \in I \cap J$.  If $k \in J$, then add $\sgn(b_k)$ since we need $+1$ if $b_k > 0$ and $-1$ if $b_k < 0$.  The indexing sets $I$ and $J$ correspond to the coefficients to which the moves (M) and (M') are applied.  Observe that $k \not \in I$ since we cannot apply (M) to the first coefficient of a continued fraction expansion.

Since $\sigma(\bit{x}_m) - \sigma(\bit{x}(0)) = \begin{cases} 0 \\ \sgn(b_k) \end{cases}$ and
\[ \sigma(\bit{x}_s(0)) - \sigma(\bit{x}(0))= \sigma(\bit{b}_1') - \sigma(\bit{x}(0)) = 1 - -2 = +3, \]
then
{\small
\begin{align*}
\left. \begin{array}{c} 0\\\sgn(b_k)\end{array} \right\}
&= \sigma(\bit{x}_m) - \sigma(\bit{x}(0))       \\
&= \sigma(\bit{x}_m) - \sigma(\bit{b}_1') + \sigma(\bit{b}_1') - \sigma(\bit{x}(N))     \\
&= 
-3 \left(\sum_{i \in I} -b_i + 
\begin{cases}
1 &\mbox{ if } 1 \in I\\
0 &\mbox{ otherwise}
\end{cases}                         
+\sum_{j \in J} b_j
\right)                                
+ 
\begin{cases}
\sgn(b_k) &\mbox{ if } k \in J\\
0 &\mbox{ otherwise}
\end{cases}
+ 3                                  \\
&= 
-3 \left(\sum_{i \in I} -b_i  + 
\begin{cases}
1 &\mbox{ if } 1 \in I \\
0 &\mbox{ otherwise}
\end{cases}
+ \sum_{j \in J} b_j
-1
\right)                           
+     
\begin{cases}
\sgn(b_k) &\mbox{ if } k \in J\\
0 &\mbox{ otherwise}
\end{cases}                         \\
&=
-3 
\left(
\sum_{i \in I} -b_i  + 
 \sum_{j \in J} b_j +
\begin{cases}
0 &\mbox{ if } 1 \in I \\
-1 &\mbox{ otherwise}
\end{cases}
\right)                            
+     
\begin{cases}
\sgn(b_k) &\mbox{ if } k \in J\\
0 &\mbox{ otherwise}
\end{cases}                         
\end{align*}
}
If $k \in J$ then 
\[
\left. \begin{array}{c} 0\\\sgn(b_k)\end{array} \right\} =
-3 
\left(
\sum_{i \in I} -b_i  + 
 \sum_{j \in J} b_j +
\begin{cases}
0 &\mbox{ if } 1 \in I \\
-1 &\mbox{ otherwise}
\end{cases}
\right)
+ \sgn(b_k)
\]
which implies we must take $\sigma(\bit{x}_m) - \sigma(\bit{x}(0)) = \sgn(b_k)$.  Thus if $k \in J$, then $N' = N - \sgn(b_k) = -\sgn(b_k)$.  

If $k \not \in J$, then similarly we must take 
which implies we must take $\sigma(\bit{x}_m) - \sigma(\bit{x}(0)) = 0$.  Thus if $k \not \in J$, then $N' = N = 0$.

In either situation,
\[ 0 =
-3 
\left(
\sum_{i \in I} -b_i  + 
 \sum_{j \in J} b_j +
\begin{cases}
0 &\mbox{ if } 1 \in I \\
-1 &\mbox{ otherwise}
\end{cases}
\right),
\]
and hence
\[ 0 =
\sum_{i \in I} -b_i  + 
 \sum_{j \in J} b_j  +
\begin{cases}
0 &\mbox{ if } 1 \in I \\
-1 &\mbox{ otherwise.}
\end{cases}
\]

\noindent{\em Subcase $N =\pm 1$.}

 Assume $[\bit{x}_m]$ or $[\bit{x}_m]$ is a MCFE for $p(\pm1)/q(\pm1)$.  Then it is obtained from $[\bit{x}_s(\pm1)] = [\bit{b}_1'']$ by applying moves (M) and (M') to nonadjacent coefficients of $[\bit{x}_s(0)]$. (Cf.\ \S\ref{subsubsec:rcfe}.)  Again, if a coefficient of $[\bit{b}_1'']$ is $\pm1$, then move (M) or (M') must be applied to it or an adjacent coefficient.

Then the difference of coefficient sums is
{\footnotesize
\begin{align*}
   \sigma(\bit{x}_m) &- \sigma(\bit{b}_1'') \\&= 
   -3(\pm 1) +
   \sum_{i \in I} -3(-b_i) +
   \begin{cases}
   -3 &\mbox{ if } 1 \in I \\
   0 & \mbox{ otherwise}
   \end{cases}          
+
   \sum_{j \in J} -3(b_j)                      
+
   \begin{cases}
   \sgn(b_k) &\mbox{ if } k \in J\\
   0 &\mbox{ otherwise}
   \end{cases}                                  \\
&=
   -3 \left(
   \pm 1  +
   \sum_{i \in I} -b_i + 
   \begin{cases}
   1 & \mbox{ if } 1 \in I \\
   0 & \mbox{ otherwise}
   \end{cases}
+
   \sum_{j \in J} b_j
   \right)                                      
 + 
   \begin{cases}
   \sgn(b_k) &\mbox{ if } k \in J\\
   0 &\mbox{ otherwise}
   \end{cases}  
\end{align*}
}
where $I \subset \{1, \dots, k-2, k\}$ and $J \subset \{1, \dots, k\}$ each not containing consecutive integers, $k \in I$, and $1 \not \in I \cap J$.   Note that we must apply move (M) to the second coefficient $-b_k \pm 1$ in $[\bit{b}_1'']$ to account for the leading coefficient of $\pm 1$, hence $k \in I$.  Also if $k \in J$, then we take $\sgn(b_k)$ since we need $+1$ if $b_k > 0$ and $-1$ if $b_k < 0$.  As noted in the previous case, the indexing sets $I$ and $J$ correspond to the coefficients to which the moves (M) and (M') are applied.

Since $\sigma(\bit{x}_m) - \sigma(\bit{x}(\pm1)) =  \begin{cases} 0 \\ \sgn(b_k) \end{cases}$ and
\begin{align*}
\sigma(\bit{x}_s(\pm1)) - \sigma(\bit{x}(\pm1)) &= \sigma(\bit{b}_1'') - \sigma(\bit{x}(\pm1)) \\ 
&= (\pm2 +1) - (\mp1 + -2) = \pm3 +3 
\end{align*}
then 
{\footnotesize
\begin{align*}
\left. \begin{array}{c} 0\\\sgn(b_k)\end{array} \right\}
&= \sigma(\bit{x}_m) - \sigma(\bit{x}(\pm1)) \\
&= \sigma(\bit{x}_m) - \sigma(\bit{b}_1'') + \sigma(\bit{b}_1'') - \sigma(\bit{x}(\pm1)) \\
&= 
-3 \left(
 \pm 1 +
\sum_{i \in I} -b_i + 
\begin{cases}
1 & \mbox{ if } 1 \in I \\
0 & \mbox{ otherwise}
\end{cases} + 
\sum_{j \in J} b_j 
\right)                                       
 +
\begin{cases}
-\sgn(b_k) &\mbox{ if } k \in J\\
0 &\mbox{ otherwise}
\end{cases} \pm3 +3                       \\
&=
-3 \left(
\pm 1 +
\sum_{i \in I} -b_i + 
\begin{cases}
  1 & \mbox{ if } 1 \in I\\
  0 & \mbox{ otherwise}
\end{cases} + 
\sum_{j \in J} b_j
\mp1 -1
\right)                                         
 +
\begin{cases}
  \sgn(b_k) &\mbox{ if } k \in J\\
  0 &\mbox{ otherwise}
\end{cases}                             \\
&=
-3 \left(
\sum_{i \in I} -b_i + 
\sum_{j \in J} b_j +
\begin{cases}
0 & \mbox{ if } 1 \in I \\
-1 & \mbox{ otherwise}
\end{cases}
\right)                               
 +
\begin{cases}
\sgn(b_k) &\mbox{ if } k \in J\\
0 &\mbox{ otherwise.}
\end{cases}
\end{align*}
}
If $k \in J$ then
\[
\left. \begin{array}{c} 0\\\sgn(b_k)\end{array} \right\} =
-3 \left(
\sum_{i \in I} -b_i + 
\sum_{j \in J} b_j +
\begin{cases}
0 & \mbox{ if } 1 \in I\\
-1 & \mbox{ otherwise}
\end{cases} 
\right)
+ \sgn(b_k)
\]
which implies we must take $\sigma(\bit{x}_m) - \sigma(\bit{x}(0)) = \sgn(b_k)$.  Thus if $k \in J$, then $N' = N - \sgn(b_k) = \pm1-\sgn(b_k) = 0$ since $\sgn(b_k) = N = \pm1$.  

If $k \not \in J$, then similarly we must take 
which implies we must take $\sigma(\bit{x}_m) - \sigma(\bit{x}(0)) = 0$.  Thus $N' = N = \pm1$.

In either situation,
\[
0 =
-3 \left(
\sum_{i \in I} -b_i + 
\sum_{j \in J} b_j +
\begin{cases}
0 &\mbox{ if } 1 \in I \\
-1 &\mbox{ otherwise}
\end{cases} 
\right)
\]
and hence
\[
0 =
\sum_{i \in I} -b_i + 
\sum_{j \in J} b_j +
\begin{cases}
0 &\mbox{ if } 1 \in I \\
-1 &\mbox{ otherwise.}
\end{cases}
\]

The only difference between the outcomes of these two subcases is in the indexing set $I$.  If $N = 0$ then $I \subset \{1, \dots, k-1\}$ containing no consecutive integers.  If $N = \pm 1$ then $I \subset \{1, \dots, k-2, k\}$ containing $k$ and no consecutive integers.  This may be reinterpreted as saying $N = 0$ if $k \not \in I$ and $N = \pm 1 = \sgn(b_k)$ if $k \in I$.  These two conditions can be consolidated by taking $I \subset \{1, \dots, k\}$ containing no consecutive integers.   Therefore a solution to the equation
\[ \tag{*}\label{*}
0 =
\sum_{i \in I} -b_i + 
\sum_{j \in J} b_j +
\begin{cases}
0 &\mbox{ if } 1 \in I \\
-1 &\mbox{ otherwise}
\end{cases}
\]
where $I$ and $J$ are subsets of $\{1, \dots, k\}$ each not containing consecutive integers and $1 \not \in I \cap J$ corresponds to an $N = 0$ or $\pm 1$, a continued fraction $[\bit{x}_m] = p(N)/q(N)$, and an $N'$ such that both $[\bit{x}_m, N'] = r(N)/s(N)$ and $\sigma(\bit{x}_m) + N' = \sigma(\bit{x}(N)) + N$.  The continued fraction $[\bit{x}_m]$ will be minimal only if moves (M) or (M') are applied to coefficients $\pm 1$ of $[\bit{x}_s(N)]$ or their neighbors.

A priori, such an $[\bit{x}_m]$ might not be of odd length.  Recall if $k \not \in I$, then $N=0$ and $[\bit{x}_s(N)]=[\bit{b}_1']$, and if $k \in I$, then $N=\sgn(b_k)$ and $[\bit{x}_s(N)]=[\bit{b}_1'']$.  Hence, similar to the method for calculating coefficient sums, the length of $[\bit{x}_m]$ may be determined from
\begin{multline*}
\len(\bit{x}_m)-\len(\bit{x}_s(N)) = 
\sum_{i \in I} (|-b_i| -2) + 
\sum_{j \in J} (|b_j| -2) \\
\shoveright{+
\begin{cases}
-\sgn(b_1) &\mbox{ if } 1 \in I \\
0 &\mbox{ otherwise} 
\end{cases}  +
\begin{cases}
-N\sgn(b_k)&\mbox{ if } k \in I \\
0&\mbox{ otherwise}
\end{cases} } \\
\end{multline*}
and
\[
\len(\bit{x}_s(N)) = 
\begin{cases}
2k+1 &\mbox{ if } k \in I \\
2k  &\mbox{ if } k \not \in I. 
\end{cases} 
\]
Since $k\in I$ implies $N=\sgn(b_k)$, $-N\sgn(b_k) = -1$.  Thus
\[ \len(\bit{x}_m) = 
\sum_{i \in I} (|-b_i| -2) + 
\sum_{j \in J} (|b_j| -2) +
\begin{cases}
-\sgn(b_1) &\mbox{ if } 1 \in I \\
0 &\mbox{ otherwise} 
\end{cases}  
+ 2k. \]
Viewing this mod $2$, $[\bit{x}_m]$ has odd length if and only if 
\[
1 \equiv
\sum_{i \in I} -b_i + 
\sum_{j \in J} b_j +
\begin{cases}
1 &\mbox{ if } 1 \in I \\
0 &\mbox{ otherwise}
\end{cases} 
+ 0 \mod 2.
\]
Therefore if and only if
\[ 0 \equiv
\sum_{i \in I} -b_i + 
\sum_{j \in J} b_j +
\begin{cases}
0 &\mbox{ if } 1 \in I \\
-1 &\mbox{ otherwise}
\end{cases} \mod 2.
\]
Hence the odd length of $[\bit{x}_m]$ is implied by the solution to equation (\ref{*}).  The conclusion of the theorem in this case follows.

\noindent{\em Cases (2), (3), and (4).}
These last three cases follow much the same via the appropriate substitutions of $[\dots,-b_2, -b_1 +1, b_1, b_2, \dots]$ by $[\dots, -b_3, -b_2, b_2 +1, b_3, \dots]$, $[\dots, -b_4, -b_3, b_3 +1, b_4, \dots]$, and $[\dots, -b_3, -b_2 +1, b_2, b_3, \dots]$ respectively.

The conclusion of Theorem~\ref{thm:generalapplication} then follows.
\end{proof}

\begin{remark}\label{rem:genuscalc}
If $I$ and $J$ are sets giving a solution to the equation of Theorem~\ref{thm:generalapplication} (1), then if $[\bit{x}_m]$ is the corresponding MCFE, the corresponding closed essential surface is $\widehat{C} = \widehat{C}(0; \bit{x}_m, N')$ where 
\[N' = \begin{cases} \sgn(b_k) & \mbox{ if } k \in I-J \\
         -\sgn(b_k) & \mbox{ if } k \in J-I \\
         0 & \mbox{ otherwise.} \end{cases}\]
The genus of $\widehat{C}$ may then be obtained from the length of $[\bit{x}_m]$:
\[g(\widehat{C}) = \frac{1}{2}(\len(\bit{x}_m)+1)-1.\]
The equation for $\len(\bit{x}_m)$ is in the proof above.

Similar formulae for the genera of surfaces corresponding to solutions of the other three equations in Theorem~\ref{thm:generalapplication} may be obtained.
\end{remark}


\section{Applications}

For the ensuing applications we recall two results from \cite{baker:sdavobkI}.  We defined a family of links $L(2n+1,\eta) \subset S^3$ such that any $L \in \K_\eta$ may be obtained by surgery on $L(2n+1, \eta)$ for some $n \in \N$ \cite[Proposition~3.1]{baker:sdavobkI}.  A continued fraction of length $n$ for the slope of $L$ confers a surgery description for $L$ on the link $L(2n+1, \eta)$.
Because the link $L(2n+1, \eta)$ is hyperbolic for $n \geq 2$, \cite[Theorem~4.1]{baker:sdavobkI} shows that knots in $\K_\eta$ with continued fractions of long length and large coefficients have large hyperbolic volume due to Thurston's Hyperbolic Dehn Surgery Theorem \cite{thurston:gt3m}. 

First, however, let us compute an example that will be of use in Corollary~\ref{cor:largevolumemanysurfaces}.
\begin{example}\label{exa:genus2again}
Let $L \in \K_\eta$ be a knot with slope $[y, -z, y]$ where $y$ and $z$ are integers $\geq 2$.  Since $[y,-z,y]$ is a SCFE with first coefficient $\neq 0$ or $1$, Theorem~\ref{thm:generalapplication}~(1) applies.  The sets $I = \{1\}$ and $J = \{3\}$ satisfy the hypotheses and provide a solution to the equation since
\[ 0 = -y + y + 0.\]
Thus there is a closed essential surface in the complement of $L$.  One may check that this is the only one.  Following Remark~\ref{rem:genuscalc}, this surface has genus $y$.   

Since $\frac{1+2z}{4z} = [2, -(z+1), 2]$, this example generalizes Example~\ref{exa:genus2sfces}.
\end{example}


\subsection{Genera of surfaces and volumes}
\begin{lemma}
If $L \in \K_\eta$ is a hyperbolic knot with $\vol(L) > \vol(L(2n+1,\eta)) > (2n+1) \cdot v$, then $S^3-N(L)$ contains no closed essential surfaces of genus less than $\frac{1}{2} (n - 1)$.
\end{lemma}

\begin{proof}
Choose $n \in \N$.  
Assume $L \in \K_\eta$ has volume greater than $\vol(L(2n+1),\eta)$.  Then since volume decreases under surgery, $L$ cannot be written as surgery on the link $L(2n+1, \eta)$ (cf.~\cite[Theorem~4.1]{baker:sdavobkI}).  Thus if $L$ has slope $\frac{x}{y}$, then any continued fraction expansion for $\frac{x}{y}$ must have length greater than $n$. 

Let $[\bit{b}] = [b_1, b_2, \dots, b_k]$ be the SCFE for $\frac{x}{y}$.  Thus $k > n$. Construct $[\bit{x}_s(N)]$ as in Lemma~\ref{lem:SCFEforp/q}.

Since $[\bit{b}]$ has length $k$, the length of $[\bit{x}_s(N)]$ must be at least $2k - 2$, the length of $[\bit{b}_4']$.  

Since the genus of a closed essential surface is related to the length of the associated MCFE of $p(N)/q(N)$ (cf.\ Remark~\ref{rem:genuscalc}), we need to see how much shorter a MCFE can be that the SCFE from which it is obtained.

Recall that any MCFE may be obtained from the corresponding SCFE by the moves (M) and (M') on non-adjacent coefficients of the SCFE.  Since length only decreases under these moves applied to a coefficient of $\pm1$, each of which decreases the length by $1$, the length of a MCFE must be at least half the length of its corresponding SCFE.  Generically, however, the length will increase.

Therefore the length of a MCFE for $p(N)/q(N)$ must be at least $\frac{1}{2}(2k-2)$ and hence greater than $k-1$.  Since $k > n$, the length of a MCFE for $p(N)/q(N)$ is greater than $n-1$.  If a MCFE $[\bit{x}_m(N)]$ for $p(N)/q(N)$ of length $k'$ corresponds to a twisted surface with meridional boundary, then the genus of the closed essential surface is $\frac{1}{2} (k'+1) - 1$.  Since $k' \geq k$, the genus of a closed essential surface in the complement of $L$ must be greater than $\frac{1}{2} (n-1) -1$.
\end{proof}

\begin{lemma}\label{lemma:manysurfaces}
There exist knots in $\K_\eta$ with arbitrarily many distinct closed essential surfaces in their complements.
\end{lemma}

\begin{proof}
Let $[\bit{b}] = [b_1, \dots, b_k]$ be the SCFE of the slope of a knot $K_0 \in \K_\eta$ with a closed essential surface in its complement such that $b_1 \neq 0$ or $1$.  By Theorem~\ref{thm:generalapplication} there are subsets $I_0$ and $J_0$ of $\{1, \dots, k\}$ such that
\[ 0 = \sum_{i \in I_0} -b_i + \sum_{j \in J_0} b_j +
\begin{cases} 0 & \mbox{ if } 1 \in I \\ -1 & \mbox{ otherwise}. \end{cases} \]
For each integer $r \geq 2$, let $[\bit{a}_r] = [a_{k+1}, a_{k+2}, \dots, a_{k+r}]$ be a SCFE and a MCFE such that $\sgn(a_{k+1}) = -\sgn(b_k)$.  Therefore $[\bit{b}, \bit{a}_r]$ is a SCFE. Let $K_r$ be the knot in $\K_\eta$ with the slope expressed by this continued fraction.

Let $\mathcal{R}_r$ be the collection of all subsets of $\{k+2, k+3, \dots, k+r\}$ that do not contain consecutive integers.  Then for every $R \in \mathcal{R}_r$,
\begin{align*}
 0
 &= \sum_{i \in I_0 \cup R} -b_i + \sum_{j \in J_0 \cup R} b_j +
\begin{cases} 0 & \mbox{ if } 1 \in I \\ -1 & \mbox{ otherwise}. \end{cases}
\end{align*}
Hence for every $R \in \mathcal{R}_r$ there is a closed essential surface in the complement of $K_r$.

Given $R, R' \in \mathcal{R}_r$ such that $R \subset R'$, then the surfaces obtained by this construction corresponding to $R$ and $R'$ have genera $g$ and $g'$ respectively such that $g < g'$ (cf.\ Remark~\ref{rem:genuscalc}).  Therefore the surfaces are distinct.  

As $r$ increases, so does the cardinality of $\mathcal{R}_r$ and size of maximal sets $R \in \mathcal{R}_r$.  Therefore as $r$ increases the number of distinct closed essential surfaces in the complement of $K_r$ increases.
\end{proof}

\begin{cor}\label{cor:largevolumemanysurfaces}
There exists arbitrarily large volume knots in $\K_\eta$ with arbitrarily many distinct closed essential surfaces in their complements.
\end{cor}

\begin{proof}
In Example~\ref{exa:genus2again}, we may freely choose positive integers $y$ and $z$ as large as we like and still obtain a knot in $\K_\eta$ of slope $[y,-z,y]$ whose complement contains an essential surface.  

Let $y=z$ and $[\bit{b}] = [z,-z, z]$.  Following the above proof of Lemma~\ref{lemma:manysurfaces} with $[\bit{a}_r] = [-z, z, -z, \dots, \pm z]$ of length $r$, the knot $K_r(z)$ of slope $[\bit{b}, \bit{a}_r]=[z, -z, z, -z, \dots, \pm z]$ contains an increasing number of distinct closed essential surfaces in its complement as $r$ increases.  Furthermore, for each $r$ as $z \rightarrow \infty$, $\vol(K_r(z)) \rightarrow \vol(L(2(r+3)+1, \eta))$.  Hence for any real numbers $V$ and $N$, there exist sufficiently large $r$ and $z$, such that $K_r(z)$ has volume greater than $V$ and more than $N$ distinct closed essential surfaces in its complement.
\end{proof}


\subsection{Small knots and volumes}

\begin{lemma} \label{lem:smallknots}
Let $\Phi_0 \geq 5$ be an integer.  Let $[\bit{b}] = [b_1, b_2, \dots, b_k]$ be a SCFE such that $b_i \equiv 0 \mod \Phi_0$ for $i = 2, \dots, k$ and $b_1 \not \equiv 0$ or $ 1 \mod \Phi_0$.  Then the knot $L \in K_{\eta}$ with slope $\frac{x}{y} = [\bit{b}]$ is small.
\end{lemma}

\begin{proof}
Because $b_1 \neq 0$ or $1$, we may apply case (1) of Theorem~\ref{thm:generalapplication}.  An essential surface in the complement of $L$ will correspond to a solution to 
\[ 0 =
   \sum_{i \in I} -b_i + 
   \sum_{j \in J} b_j +
   \begin{cases}
      0 &\mbox{ if } 1 \in I \\
      -1 &\mbox{ otherwise}
   \end{cases}
\]
subject to various constraints on the sets $I$ and $J$.  A solution to this equation will hold true $\mod \Phi$.  By choosing $b_i \equiv 0 \mod \Phi$ for $i \geq 2$, this equation becomes
\[
 0 \equiv   
   \begin{cases}
      -b_1 &\mbox{ if } 1 \in I\\
      b_1 -1 &\mbox{ if } 1 \in J\\
      -1 &\mbox{ otherwise}
   \end{cases} \mod \Phi.
\]
This equation can only hold true if $b_1 \equiv 0$ or $1 \mod \Phi$.  Having chosen $b_1$ otherwise, this equation and hence the original equations have no solution.  Therefore there can be no closed essential surfaces in the complement of $L$.
\end{proof}

\begin{cor}\label{cor:largevolumenosurfaces}
There exists arbitrarily large volume knots in $\K_\eta$ with no closed essential surfaces in their complements.
\end{cor}

\begin{proof}
For each integer $n \geq 2$ choose an integer $\Phi_n >> 5$ that is ``sufficiently large'' as at the end of the proof of \cite[Theorem~4.1]{baker:sdavobkI}.  Consider the sequence of knots $K_n \in \K_{\eta}$ represented by the simple continued fractions of length $n$ 
\[[\bit{k}_n] = [\Phi_n +2, -\Phi_n, +\Phi_n, -\Phi_n, +\Phi_n, \dots, \pm\Phi_n].\]
By \cite[Theorem~4.1]{baker:sdavobkI}, $\vol(K_n) \to \infty$ as $n \to \infty$.  Thus the set of volumes $\{\vol(K_n)\}$ is unbounded.

 All the coefficients of $[\bit{k}_n]$ except the first are congruent to $0 \mod \Phi_n$.  The first coefficient is $\Phi_n +2 \equiv 2 \mod \Phi_n$ and is not congruent to $0$ or $1$.  Since the coefficients of $[\bit{k}_n]$ alternate sign, it is a simple continued fraction.  By Lemma~\ref{lem:smallknots}, the knots $K_n$ are all small.
\end{proof}

\begin{remark}
It is well known that every two-bridge knot is small and that there are hyperbolic two-bridge knots of arbitrarily large volume.  No hyperbolic two-bridge knot, however, admits a lens space surgery \cite{takahashi}.
\end{remark}

\bibliography{MathBiblio}
\bibliographystyle{plain}
\end{document}